\documentclass[12pt]{article}

\usepackage[final]{epsfig}
\usepackage{graphics}
\usepackage{epsfig,color}
\usepackage{amsmath}
\usepackage{amsfonts}
\usepackage{latexsym}
\usepackage{amssymb}

\newtheorem{prop}{Proposition}
\newtheorem{defin}{Definition}
\newtheorem{thm}{Theorem}
\newtheorem{corol}{Corollary}
\newtheorem{lem}{Lemma}

\newtheorem{rem}{Remark}

\def\qed{\hfill\vrule height 5pt width 5 pt depth 0pt\\}

\newcommand{\A}{{\cal A}}
\newcommand{\BB}{{\cal B}}

\newcommand{\HHH}{{\bf H}}
\newcommand{\HT}{{\bf H}^2}

\newcommand{\LL}{{\cal L}}

\newcommand{\N}{{\bf N}}

\newcommand{\PP}{{\cal P}}

\newcommand{\QQ}{{\cal Q}}
\newcommand{\R}{\bf R}
\newcommand{\RR}{\bf R^2}

\newcommand{\Z}{{\bf Z}}

\newcommand{\ar}{\mbox{area}}
\newcommand{\bo}{\partial}
\newcommand{\ess}{\mbox{ess}}
\newcommand{\inter}{\mbox{int}}
\newcommand{\per}{\mbox{perimeter}}
\newcommand{\sph}{\bf S^2}

\begin{document}

\bibliographystyle{plain}

\title
{Complexity of piecewise convex transformations in two
dimensions, with applications to polygonal billiards}
\author{Eugene Gutkin
\thanks{UCLA, Math. Dpt., Los Angeles, CA 90095, USA  and IMPA,
Estrada Dona Castorina 110, 22460-320 Rio de Janeiro RJ, Brasil.
{\em E-mail}: egutkin@math.ucla.edu,gutkin@impa.br}
\ and Serge Tabachnikov
\thanks{Department of Mathematics, Pennsylvania State University,
University Park,
PA 16802. {\em E-mail}: tabachni@math.psu.edu}}
\maketitle

\begin{abstract}                             
ABSTRACT
\noindent We introduce the class of piecewise convex transformations,
and study their complexity. We apply the results to the complexity
of polygonal billiards on surfaces of constant curvature.
\end{abstract}

\vspace{5mm}

\section*{Introduction}         \label{intro}
The following situation frequently occurs in geometric dynamics.
There is a phase space $X$, a transformation $T:X\to X$; and there is
a finite decomposition  $\PP:X=X(a)\cup X(b)\cup\cdots$.
Let $\A=\{a,b,\dots\}$ be the corresponding alphabet.
A phase point $x\in X$ is regular if every element of the orbit
$x,Tx,T^2x,\dots$ belongs to a unique atom of $\PP$.
Suppose that $x\in X(a),Tx\in X(b)$, etc. The corresponding
word $a\,b\,\cdots$ is the {\em code} of $x$.
Let $\Sigma(n)$ be the set of words in $\A$ of length $n$ obtained
by coding points in $X$. The positive function $f(n)=|\Sigma(n)|$ is the
associated {\em complexity}. Its behavior as $n\to\infty$\footnote{
The function $f(\cdot)$ can be bounded, or grow polynomially,
or grow exponentially, etc.} is an important characteristic of
the dynamical system in question. The following examples have motivated
our study.

\medskip

\noindent{\bf Example A}.  Let $P\subset\RR$ be a polygon
with sides $a,b,\dots$,  and let $X$ be the
phase space of the {\em billiard map} $T_b$ in $P$.
The coding generated by the corresponding decomposition
$\PP:X=X(a)\cup X(b)\cup\cdots$ is the traditional coding
of billiard orbits by the sides they hit \cite{Ta95}.
Basic questions about its complexity are open \cite{Gu03}.

\vspace{3mm}

\noindent{\bf Example B}. Let $P\subset \RR$ be a convex polygon
with vertices $a,b,\dots$.
The complement $X= \RR-P$ is the phase space of the
{\em outer billiard} $T_o$ about $P$. (It is also
called the {\em dual billiard}. See \cite{Ta95}).
The conical regions bounded the singular lines of $T_o$
form the natural decomposition $\PP:X=X(a)\cup X(b)\cup\cdots$.
In $X(a),X(b),\dots$ the mapping $T_o$ is the symmetry about
$a,b,\dots$. The decomposition  $\PP$ yields
the coding of outer billiard orbits by the vertices
they hit.

\medskip

We will study the complexity of ($2$-dimensional)
{\em piecewise convex transformations}.
This is a wide class of geometric dynamical systems; it contains the  
examples
above.
Our setting is as follows. (For  simplicity of exposition,
we restrict our attention to two dimensions.)

Let $X$ be a {\em geodesic surface}, and let $\Gamma\subset X$ be a
finite {\em geodesic graph}. A subset $Y\subset X$ is {\em convex}
if for any $x,y\in Y$ there is a unique geodesic in $Y$ with endpoints
$x,y\in Y$. Suppose that the closed faces  $X(a),\, X(b),\,,\dots$
of $\Gamma$ are convex, and let $\PP:X=X(a)\cup X(b)\cup\cdots$
be the corresponding decomposition. We say that $\PP$ is a
{\em convex partition} of $X$. A differentiable mapping
$T:Y\to X$ is convex if it sends geodesics to geodesics.
Suppose that $T:X\to X$ is a convex diffeomorphism
on $X(a),\, X(b),\,,\dots$, and that $T(X(a)),\, T(X(b)),\,,\dots$
also form a convex partition of $X$. We say that the triple
$(X,T,\PP)$ is a {\em piecewise convex transformation} with
the {\em defining partition} $\PP$.

Making further assumptions on $(X,T,\PP)$, we obtain
more specialized classes of transformations, e.g.,
piecewise isometries, piecewise affine mappings, etc.
The case when the faces
of $\Gamma$ are convex euclidean polyhedra, and $T$ is
isometric on them
arises in task scheduling problems \cite{Ad-etal04}.

In Section~\ref{basic} we develop geometric and combinatorial techniques
to study the complexity of piecewise convex transformations.
In the rest of the paper we apply these results to the inner and
outer polygonal billiards on (simply connected)
surfaces of constant curvature $\chi=0,1,-1$.
One of the goals of this work is to develop a uniform approach to these
dynamical systems. While there is a vast literature on the parabolic  
case
($\chi=0$), the elliptic ($\chi=1$) and the hyperbolic ($\chi=-1$)
cases have been studied only sporadically.

Let $M$ be a surface of constant curvature, and let $P\subset M$ be
a polygon. In Section~\ref{innerbill} we cast the (inner) billiard map  
in
$P$ as a piecewise convex transformation $(X,T,\PP)$. See
Theorem~\ref{bil_piece_prop}.
We do it simultaneously for all curvatures,
and without making additional assumptions on $P$.
Thus, we don't assume that $P$ is convex or simple, and have to pay a  
price
for this. The partition $\PP$ is finer than
the ``natural" one; the coding it generates is more refined  than the  
standard coding by sides \cite{Ta95}.
We develop a dictionary between the language of
piecewise convex transformations and that of billiard orbits.
In the remaining part of the paper we come back to
Examples A, B and study the complexity of natural coding.

In Section~\ref{cases} we investigate
the complexity of inner billiard orbits in a polygon $P$ on
a surface of constant curvature $\chi$.
In section~\ref{arbit_sub} $P,\chi$ are arbitrary; later on
we specialize to convex $P$ but arbitrary $\chi$,
to $\chi=0$, and to $\chi=\pm 1$ respectively. Below we formulate the  
main results.

\medskip

\noindent The side complexity of billiard orbits in any rational
euclidean polygon grows at most cubically; see
Theorem~\ref{rat_complex_thm}.

\noindent The side complexity of billiard orbits in any spherical
polygon grows subexponentially; see
Theorem~\ref{spheresubexp}.

\noindent The side complexity of billiard orbits in any hyperbolic
(i.e., $\chi=-1$) polygon grows exponentially; the exponent in question  
is
the topological entropy of the billiard map; see
Theorem~\ref{hyperb_exp}.

Section~\ref{outerbill} is the outer billiard counterpart of
Section~\ref{cases}. Here are some of its results.

\medskip

\noindent For $\chi=0$ and arbitrary polygon (resp. rational polygon)
we obtain polynomial bounds from above and below (resp. quadratic  
asymptotics)
for the complexity; see Theorem~\ref{upperbound} and  
Theorem~\ref{rational}
respectively.

\noindent For $\chi=1$ complexity grows subexponentially;
see Theorem~\ref{spheresubexp1}.

\noindent For $\chi=-1$  and arbitrary polygon, we obtain linear lower  
bound for
complexity,
which is sharp: for the so-called large polygons complexity grows  
linearly;
see Theorem~\ref{hyplowbd}.

\vspace{5mm}

\noindent {\bf Notes and references.} We are grateful to the Research in Pairs  
program in Oberwolfach and to the Shapiro visiting program at Penn State for
their support. The second author was partially supported by an NSF  
grant. The work \cite{CHT03,Be03} is the predecessor of our
Section 3.\footnote{See especially Remark 7 in \cite{Be03}.}
For additional results on polygonal billiards see  
\cite{GS92,Ka87,Tr98,Sh-Vi}.

\section{Piecewise convex transformations}   \label{basic}
\subsection{Convex geodesic surfaces}   \label{basic_sub}
We will introduce a class of maps of {\em geodesic spaces}.
In order to simplify our exposition, we will restrict it to two
dimensions, i.e., to  {\em geodesic surfaces}. Let  $M$ be one.
Then $M$ is a topological surface.
It may have a boundary, $\bo M$, and a finite number
of {\em cone points}.  The surface  $M$ is endowed
with a collection of {\em geodesics}, satisfying the standard
properties. In particular, $\bo M$ is a finite union of geodesics.
Let $x,y\in M$ be any pair of points.  Then
$M$ is {\em geodesically convex} if there is a {\em unique
shortest geodesic} $\gamma=\gamma(x,y)\subset M$ joining
them.\footnote{Although the uniqueness condition may seem too  
restrictive,
it is crucial for our study of iterations of piecewise convex  
transformations.
See section~\ref{appli_sub} below.}

\vspace{3mm}

\noindent {\bf Examples.} i) The set $\RR$ endowed with a projectively
flat  (e.g., Minkowski) metric is a geodesic surface.
Its geodesics are straight lines. Let
$M\subset\RR$ be a closed, bounded convex region.
Then it is a geodesic surface
iff $M$ is a polygon. The geodesic
$\gamma(x,y)$ is the segment with endpoints $x,y$.\\
ii) Let $\sph$ be the round sphere, and let $M\subset\sph$
be a convex spherical polygon. Then $M$ is a geodesic surface
iff it does not contain antipodal points. This holds iff
$M$ is contained in an open hemisphere.

\vspace{3mm}

In what follows we will often say ``convex geodesic surface"
instead of ``geodesically convex geodesic surface".
Let $X$ be a geodesic surface. A finite graph $\Gamma$ drawn on $X$
is a {\em geodesic graph} if the edges of $\Gamma$ are geodesics.
We will assume that the cone points of $X$
are contained in the set of vertices of $\Gamma$,
and that the latter are nondegenerate. (A vertex is degenerate
if it has two adjacent edges, and they are colinear.)

Let $X(a),X(b),\dots$ be the (open) faces of $\Gamma$, and let
$P_a,P_b,\dots\subset X$ be their closures. We will say that
$\Gamma$ is a {\em convex geodesic graph} if the surfaces
$P_a,P_b,\dots\subset X$ are geodesically convex. We will refer to
the data $(X,\Gamma)$ as a {\em piecewise (geodesically) convex
geodesic surface}. We will also say that
$\PP: X=P_a\cup P_b,\dots$ is a {\em convex geodesic partition} of
the geodesic surface $X$.
The boundary $\bo\PP\subset  X$ is spanned
by the edges and vertices of $\Gamma$. We will also refer to it
as the {\em support} of
$\Gamma$ and denote by $<\Gamma>$.  There is a $1-1$
correspondence between  convex geodesic
graphs drawn on $X$ and  convex geodesic partitions of $X$. The latter
concept   provides an alternative
approach to the material below \cite{GH97} . However, the language of  
convex geodesic
graphs is more suitable for our purposes, and we will use it in what  
follows.

Let $X$ be a convex geodesic surface, and let $Y$ be an
arbitrary geodesic surface.
A diffeomorphism (not necessarily surjective) $T:X\to Y$ is a {\em
geodesically convex transformation} if it sends geodesics to geodesics.
Then $T(X)=Z\subset Y$ is a
convex geodesic surface as well, and the inverse map $T^{-1}:Z\to X$ is  
a geodesically convex transformation.

\subsection{ Piecewise convex transformations: iterations,
coding, and complexity}
\label{convex_sub1}
We will now introduce a class of dynamical systems which, on one
hand, is sufficiently general
to include several interesting examples, and on the other, is special
enough to allow a common geometric framework. Let  
$X,\Gamma,\PP=\PP(\Gamma)$
be as above. Let $p$ be the number of atoms of $\PP$ and write
$\PP:X=\cup_{i=1}^p P_i$. Suppose that for
$1\le i \le p$ there is a convex transformation $T_i:P_i\to X$. Set
$Q_i=T_i(P_i)$.  Suppose that $X=\cup_{i=1}^p Q_i$ and that the open
sets $\inter(Q_i)$ are
pairwise disjoint. Then $\QQ:X=\cup_{i=1}^pQ_i$ is a convex geodesic  
partition,
and $\bo\QQ$ is the support of a convex geodesic graph,
$\Gamma^{-1}$, drawn on $X$.

These data determine a
{\em piecewise convex transformation} $T:X\to X$ with  the {\em defining
partition} $\PP$, and we will use the notation  $(X,T,\,\PP)$ for it.
The inverse of $(X,T,\,\PP)$ is also
a  piecewise convex transformation. Its  defining partition is
$\QQ=T(\PP)$. Thus in our
notation, $(X,T,\,\PP)^{-1}=(X,T^{-1},\QQ)$.

Invertible piecewise isometric (as well as affine, or projective) maps
are  examples of piecewise  convex transformations.
Let $M$ be a simply connected surface of constant curvature, and let
$P\subset M$ be a geodesic polygon. In Sections~\ref{cases}
and~\ref{outerbill} we will put the inner and the outer billiard about
$P$ into the framework of piecewise convex transformations.

Denote by $\LL$ the {\em full shift space} on the {\em alphabet}
$\A=\{a,b,c,\dots\}$ of the faces of $\Gamma$. A point $x\in X$ is
{\em regular} if $x,Tx,T^2x,\dots,$ belong to open faces of $\Gamma$.
Let $X_{\infty}\subset X$ be the set of regular points.
Assigning to  $x\in X_{\infty}$ the sequence of faces of $\Gamma$  
containing the
consecutive elements $x,Tx,T^2x,\dots,$ we obtain the {\em coding map}
$\sigma:X_{\infty}\to\LL$. Set $\Sigma=\sigma(X_{\infty})$, and let
$\Sigma(n)$ be the set of words of length $n$ that occur in $\Sigma$.
The function $f(n)=|\Sigma(n)|$ is the {\em complexity} of $(X,T,\PP)$
{\em with respect to the defining partition}.

Let  $\Gamma',\Gamma''$ be two
geodesic graphs on $X$.  Their {\em join} $\Gamma'\vee\Gamma''$ is the
(unique) geodesic graph such that
$<\Gamma'\vee\Gamma''>\ =\ <\Gamma'>\,\cup\,<\Gamma''>.$
We write $\Gamma''\prec\Gamma'$ if $\Gamma'\vee\Gamma''=\Gamma''$.
Recall that $\PP'\vee\PP''$ denotes the join of partitions, and that
$\PP''\prec\PP'$
holds if $\PP'\vee\PP''=\PP''$.
Let $\PP'=\PP(\Gamma'),\PP''=\PP(\Gamma'')$.  If
$\Gamma',\Gamma''$ are convex  geodesic graphs,  then so is
$\Gamma'\vee\Gamma''$, and
$
\PP(\Gamma'\vee\Gamma'')\ =\ \PP(\Gamma')\,\vee\, \PP(\Gamma'').
$
Moreover, $\Gamma''\prec\Gamma'$ iff   $\PP''\prec\PP'$.

Set $\Gamma_1=\Gamma,\Gamma_2=\Gamma_1\vee T^{-1}(\Gamma),\dots,
\Gamma_{n+1}=\Gamma_n\vee T^{-n}(\Gamma),\dots$. The convex geodesic  
graphs
$\Gamma_k,\,k=1,2,\dots$ form an increasing tower with respect to the
relation $\prec$, and we set $S_k=<\Gamma_k>$, $\PP_k=\PP(\Gamma_k)$.
The {\em singular set} $S_{\infty}=\cup_{k=1}^{\infty}S_k=X\setminus  
X_{\infty}$
is a
countable (at most) union of geodesics. Note that $(X,T^n,\,\PP_n)$ is a
piecewise convex transformation for $n\ge 1$.

Let $x=x_1,Tx=x_2,\dots,T^{n-1}x=x_n$ be a {\em finite orbit} such that
the points $x_1,x_2,\dots,x_n$ belong to open faces of $\Gamma$. We say  
that
$x_1,x_2,\dots,x_n$ is a {\em regular orbit of length $n$} and
denote by $\sigma(x_1,x_2,\dots,x_n)$ the corresponding word on the  
alphabet
$\A$; it is the {\em code of $x_1,x_2,\dots,x_n$}. Then
$\Sigma(n)$ is the set of codes of regular orbits of length $n$; the
orbit $x=x_1,Tx=x_2,\dots,T^{n-1}x=x_n$ is regular iff $x$ belongs to  
an open
face of the graph $\Gamma_n$.

The proposition below summarizes the discussion.
\begin{prop}  \label{sing_strat_prop} Let $(X,T,\,\PP)$ be a piecewise  
convex
transformation, and let $<\Gamma>=\bo\PP$. Then:\\
1. There is a sequence $\Gamma_k,\,k\ge 1$ of convex geodesic graphs in  
$X$
such that
$$
\cdots\prec\Gamma_{n+1}\prec\Gamma_n\prec\cdots\prec\Gamma_1=
\Gamma;
$$
2. Set $\PP_n=\PP(\Gamma_n),\,n\ge 1$. Then $\PP_n$ is a convex geodesic
partition, and
the $n$th iteration of $(X,T,\PP)$ is a piecewise convex
transformation with the defining partition $\PP_n$, i. e.,
$(X,T,\PP)^n=(X,T^n,\,\PP_n)$.\\
3. There is a natural bijection between the set $\Sigma(n)$
and the set of atoms of the partition $\PP_n$. Thus, the complexity
$f(n)$ is the number of faces of the graph $\Gamma_n$.
\end{prop}

\subsection{A combinatorial lemma}   \label{comb-top_sub}
We will need a general proposition that concerns the combinatorics
of graphs.
Let $X$ be a piecewise geodesic surface. Let $\Gamma',\Gamma''$ be
geodesic graphs drawn on
$X$ and set $\Gamma=\Gamma'\vee\Gamma''$.  (We do not assume that the  
faces of
$\Gamma',\Gamma''$ are convex.)

\begin{figure}[htbp]
\begin{center}
\input{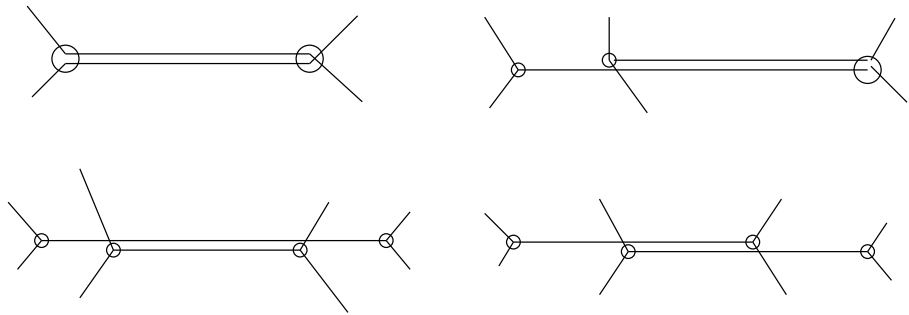}
\caption{Overlapping of edges of two graphs}
\label{fig1}
\end{center}
\end{figure}

Denote by $F',F'',E',E'',V',V''$  the respective sets of faces, edges  
and
vertices.  Let $E,F,V$ be the sets of faces, edges and vertices of  
$\Gamma$. If
$e'\in E',e''\in E''$ intersect non-transversally, then they
(partially) overlap.
There are $4$ ways in which this can happen. See figure~\ref{fig1}.  
Denote by
$c(\Gamma',\Gamma'')$ the number of the overlapping pairs of edges.

\begin{lem}  \label{comb-top_lem}
Let $\Gamma',\Gamma''$ be
geodesic graphs drawn on  a piecewise geodesic surface $X$. Denote by
$\chi$  the
Euler characteristic of $X$. Let $V'_d,V''_d$ be the sets of vertices   
of
$\Gamma',\Gamma''$ that are disjoint from the other graph, and set
$V_{\ess}=V-V'_d-V''_d$. Then
\begin{equation}  \label{comb-top_eq}
|F|-|F'|-|F''| +\chi = |V_{\ess}|  - c(\Gamma',\Gamma'').
\end{equation}
\end{lem}
{\bf Proof}. Any graph $A$ drawn on  $X$ satisfies
\begin{equation}  \label{euler_eq} |F(A)|-|E(A)|+|V(A)|=\chi(X).
\end{equation}
      Applying (\ref{euler_eq}) to graphs $\Gamma',\Gamma'',\Gamma$, we  
obtain
for the left hand side of (\ref{comb-top_eq})
\begin{equation}
\label{euler_eq1}
|F|-|F'|-|F''| +\chi =   (|E|-|E'|-|E''|)+(|V'|+|V''|-|V|).
\end{equation}

Denote by $e'_i,e''_j$ the edges of $\Gamma',\Gamma''$ respectively. Let
$a_i',a''_j$ be the number of vertices of $\Gamma'',\Gamma'$ that are
located in
the interior of $e'_i,e''_j$ respectively. Let $b_i',b_j''$ be the
number of times
that the interior of $e'_i,e_j''$ transversally intersects the
interior of an edge
of $\Gamma'',\Gamma'$ respectively. Then $e'_i,e''_j$ contribute
$a_i'+b_i'+1,a_j''+b_j''+1$  edges to $E$ respectively. Taking the
overlapping
into account, we obtain
\begin{equation}     \label{edges_eq1}
|E|=\sum_ia_i'+\sum_ib_i'+|E'|+\sum_ja_j''+ \sum_jb_j''+|E''| -
c(\Gamma',\Gamma'').
\end{equation}
Let $V_c$ be the set of common vertices of $\Gamma',\Gamma''$.  Let
$V'_e,V''_e$
be the sets of vertices  of $\Gamma',\Gamma''$ that are in the
interior of edges
of the other graph. Then
\begin{equation}     \label{vertices_eq1}
|V'|=|V'_e|+|V'_d|+|V_c|,\quad |V''|=|V''_e|+|V''_d|+|V_c| .
\end{equation}
      Besides
\begin{equation}     \label{vertices_eq2}
\sum_ia_i'=|V''_e|,\quad \sum_ja_j''=|V'_e|.
\end{equation}
      Let $V_n=V\setminus(V'\cup V'')$ be the set of ``new" vertices of
$\Gamma$. Then
\begin{equation}     \label{vertices_eq3}
\sum_ib_i'= \sum_jb_j''=|V_n|.
\end{equation}
We also have
\begin{equation}     \label{vertices_eq4}
|V|=|V'_e|+|V''_e|+|V'_d|+|V''_d|+|V_c|+|V_n|.
\end{equation}

 From equations  
 (\ref{edges_eq1}),(\ref{vertices_eq1}),(\ref{vertices_eq2}) and
(\ref{vertices_eq3}), we obtain
\begin{equation}     \label{edge_vert_eq}
|E|-|E'|-|E''|= |V'_e|+|V''_e|+2|V_n| -c(\Gamma',\Gamma'').
\end{equation}
      By equation (\ref{vertices_eq4}),
\begin{equation}     \label{vertices_eq5}
|V'|+|V''|-|V|=|V_c|-|V_n|.
\end{equation}
      Substituting (\ref{edge_vert_eq}),(\ref{vertices_eq5}) into
(\ref{euler_eq1}),
and using (\ref{vertices_eq4}) again, we obtain the claim. \qed

\vspace{3mm}

The expression for  
$|F(\Gamma'\vee\Gamma'')|-|F(\Gamma')|-|F(\Gamma'')|$ in
Lemma~\ref{comb-top_lem} involves information about
vertices of $\Gamma',\Gamma''$ which is sometimes not
available. The following corollary of Lemma~\ref{comb-top_lem} is  
useful.

\begin{prop}  \label{comb_bounds_prop}
Let $\Gamma',\Gamma''$ be geodesic graphs drawn on  a piecewise
geodesic surface $X$
of Euler characteristic  $\chi$. Let $\Gamma=\Gamma'\vee\Gamma''$, and  
let
the notation be as above. Then
\begin{equation}     \label{ineq_eq}
|V_n|   \le |F|-|F'|-|F''| + \chi + c(\Gamma',\Gamma'') \le  |V|.
\end{equation}
\end{prop}
{\bf Proof.} By Lemma~\ref{comb-top_lem} and equation
(\ref{vertices_eq4}), the
quantity to  estimate in (\ref{ineq_eq}) is equal to
$|V'_e|+|V''_e|+|V_c|+|V_n|$.
\qed

\vspace{3mm}

The following special cases of  Lemma~\ref{comb-top_lem} and
Proposition~\ref{comb_bounds_prop} will be especially useful.

\begin{figure}[htbp]
\begin{center}
\input{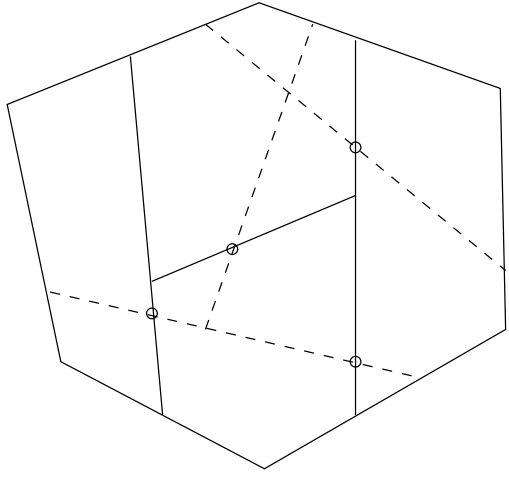}
\caption{Illustrating formula (\ref{transv_eq}): $|F|=12, |F'|=4,  
|F''|=5,
|V_{\ess}|=4$}
\label{fig1.5}
\end{center}
\end{figure}

\begin{corol}  \label{comb_bounds_cor}
Let $X$ be a geodesic surface homeomorphic to the open disc. Let
$\Gamma',\Gamma''$ be
geodesic graphs drawn on $X$ such that their edges intersect only
transversally.
Let $\Gamma=\Gamma'\vee\Gamma''$ and let the sets
$F,F',F'',V,V_n,V'_d,V''_d,V_{\ess}$ be as above. Then
\begin{equation}     \label{transv_eq}
|F|-|F'|-|F''| + 1 = |V_{\ess}|
\end{equation}
      (see figure \ref{fig1.5}), and
\begin{equation}     \label{transv_ineq}
|V_n| \le |F|-|F'|-|F''| + 1 \le |V|.
\end{equation}
      If $V_d'=V_d''=\emptyset$ then
\begin{equation}     \label{transv_eq1}
|F|-|F'|-|F''| + 1 = |V|.
\end{equation}
\end{corol}
{\bf Proof}. By our assumptions, $c(\Gamma',\Gamma'')=0$ and
$\chi=1$. Hence
Lemma~\ref{comb-top_lem} and Proposition~\ref{comb_bounds_prop} yield
(\ref{transv_eq}) and  (\ref{transv_ineq}) respectively. If
$V_d'=V_d''=\emptyset$ then $V=V_{\ess}$ and  (\ref{transv_eq}) becomes
(\ref{transv_eq1}). \qed

\subsection{Geometric formula for complexity of piecewise convex
transformations}
\label{appli_sub}
Let $\Sigma$ be a language on a finite alphabet $\A$, and let
$\Sigma(n)$ be
the set  of words of length $n$ in $\Sigma$. The {\em complexity} of
$\Sigma$ is the function
$f(n)=|\Sigma(n)|$. Set $\varphi(n)=f(n+1)-f(n)$ and
$\psi(n)=\varphi(n+1)-\varphi(n)$.
We will refer to the functions $\varphi(\cdot),\psi(\cdot)$ as the
{\em first, second
differences of complexity}, respectively.

Our approach is based on Cassaigne's formula for $\psi(\cdot)$  
\cite{Ca97}.
Denote by $m_l(w),m_r(w)$
the number of left, right one-letter extensions of the word $w$  
respectively.
We will assume, following \cite{Ca97}, that $m_l(w),m_r(w)\ge 1$
for any $w\in\Sigma$.
A word is {\em bispecial} if $m_l(w),m_r(w)>1$.
Let $m_b(w)$ be the  number of extensions of the type $awb$ where
$a,b\in \A$.
Let $\BB\subset\Sigma$ be the set of bispecial words, and set
$\BB(n)=\BB\cap\Sigma(n)$.
The {\em Cassaigne index} is defined by
\begin{equation}     \label{Cassindex}
\mu(w)=m_b(w)-m_l(w)-m_r(w)+1.
\end{equation}
Note that if $w$ is not bispecial then $\mu(w)=0$.
The Cassaigne formula says
\begin{equation}  \label{Casform_eq}
\psi(n)=\sum_{w\in\BB(n)} \mu(w) = \sum_{w\in\Sigma(n)} \mu(w).
\end{equation}
We define the cumulative index by $\mu(n)=\sum_{w\in\Sigma(n)}\mu(w)$  
for
$n\ge 1$, and $\mu(0)=0$. Set
\begin{equation}  \label{Casind_eq}
M(n)=\sum_{k\le n}\mu(k).
\end{equation}

\begin{lem}   \label{complviaCass}
The complexity of a piecewise convex transformation satisfies
\begin{equation}  \label{count_eq1}
f(n)=f(1)+(n-1)(f(2)-f(1))+\sum_{k\le n-2}M(k).
\end{equation}
\end{lem}
{\bf Proof.} Denote by $g(n)$ the right hand side of
(\ref{count_eq1}). Since, by  
(\ref{Casform_eq},\ref{Casind_eq}),
the second differences of $f(\cdot),g(\cdot)$ are equal,
$f(n)-g(n)$ is linear in $n$. But $f(1)-g(1)=f(2)-g(2)=0$. \qed

Let $(X,T,\PP)$ be a piecewise convex transformation, and let $\Sigma$
be its coding language. For $w\in\Sigma(n)$ let $X(w)\subset X$ be the
corresponding open
face of the graph $\Gamma_n$. We denote by $\Gamma(w)$ the restriction  
of
$\Gamma_{-1}\vee\Gamma_{n+1}$ to $X(w)$. Let $V_{\ess}(w)$ (resp.  
$OE(w)$)
be the set of essential vertices (resp. pairs of overlappings edges) for
$\Gamma(w)$.

\begin{lem}   \label{indexint}
For any $w\in\Sigma$ we have
\begin{equation}     \label{second_eq1}
\mu(w) = |V_{\ess}(w)|-|OE(w)|.
\end{equation}
\end{lem}
{\bf Proof.} Denote by $\Gamma',\Gamma''$ the restrictions of
$\Gamma_{n+1},\Gamma_{-1}$ to $X(w)$ respectively. Then in the notation  
of
Lemma~\ref{comb-top_lem}, $m_b(w)=|F|,
m_r(w)=|F'|,m_l(w)=|F''|$. Since $X(w)$ is contractible, $\chi=1$ in  
equation
(\ref{comb-top_eq}). Thus, the left hand side of
equation~(\ref{comb-top_eq}) is $\mu(w)$.
But its right hand side is $|V_{\ess}(w)|-|OE(w)|$. \qed

We will use the notation $c(w)=|OE(w)|$. For $n\ge 1$ set
\begin{equation}  \label{define_eq1}
V_{\ess}(n)=\cup_{w\in\Sigma(n)}V_{\ess}(w),\ v(n)=|V_{\ess}(n)|,\
V(n)=\sum_{k\le n}v(k),
\end{equation}
and
\begin{equation}  \label{vertices_eq}
OE(n)=\cup_{w\in\Sigma(n)}OE(w),\ c(n)=|OE(n)|,\ C(n)=\sum_{k\le n}c(k).
\end{equation}
Thus, $v(n)$ (resp. $c(n)$) is the number of ``new" essential vertices
(resp. ``new" edge overlappings) of the graph  
$\Gamma_{-1}\vee\Gamma_{n+1}$,
while $V(n)$ (resp. $C(n)$) is the total number of essential vertices
(resp. edge overlappings)  of $\Gamma_{-1}\vee\Gamma_{n+1}$.
Note that only bispecial words contribute to these numbers.

\begin{prop} \label{geomCassform}
Let $(X,T,\PP)$ be a piecewise convex transformation and let
$\PP_1=\PP,\PP_2,\dots$
be the corresponding sequence of convex partitions. Then the complexity  
of
$(X,T,\PP)$ satisfies
\begin{equation}  \label{complex_eq}
f(n)=|\PP_1|+(n-1)(|\PP_2|-|\PP_1|)+\sum_{k\le n-2}V(k)-\sum_{k\le  
n-2}C(k).
\end{equation}
\end{prop}
{\bf Proof.} We have $f(1)=|\PP_1|,\,f(2)=|\PP_2|$. By  
Lemma~\ref{indexint},
$\mu(n)=v(n)-c(n)$, hence $M(n)=V(n)-C(n)$. The claim now follows from
equation~(\ref{count_eq1}).\qed

\section{Billiard map as a piecewise convex transformation: the  
dictionary}
\label{innerbill}

Let $M$ be a complete simply connected surface of constant curvature,\footnote{We
normalize the curvature to $0$ or $-1$ or $1$. Thus $M$ is either the  
euclidean
plane or the hyperbolic plane or the unit sphere. We will refer to these
geometries as parabolic, hyperbolic and elliptic respectively.}
and let $P\subset M$ be a
geodesic polygon. In order to cast the billiard map in $P$ as a
piecewise convex transformation simultaneously for all three cases at  
hand,
we will  use the models of the
hyperbolic and elliptic geometries where the geodesics
are straight lines in the euclidean plane.
(We will refer to them as the {\em projective models}.)

For $\chi=-1$, this is the Klein-Beltrami  model of hyperbolic
geometry. The hyperbolic plane is represented by the open unit disc,
the unit circle $S$  is the ``circle at infinity",
and geodesics are the chords of this disc.
Let $x,y$ be distinct points in the disc, and let
$a,b$ be the intersection points of the line $xy$ with  $S$.
We identify the oriented line $xy$ with $\R$, and hence $x,y,a,b$
with real numbers. The cross-ratio of the four points
is given by
$$
[a,x,y,b]={{(y-a)(b-x)}\over{(x-a)(b-y)}}.
$$
The distance between points satisfies
$2d(x,y) = \ln [a,x,y,b].$ In this model,
isometries of the hyperbolic geometry are the
projective transformations of the euclidean plane, preserving
the unit disc.

Projective model of the elliptic geometry (i.e.,  $\chi=+1$)
is as follows. We restrict our attention to an open hemisphere;
consider for concreteness the northern hemisphere. Its central
projection to the tangent plane at $(0,0,1)$  is a
surjective diffeomorphism that takes spherical geodesics to
euclidean straight lines.

In order to endow the {\em billiard map phase space} $X=X(P)$
with the structure of a geodesic surface in a uniform fashion,
we will impose a slight restriction
on $P$. Let $P\subset M$ be a geodesic polygon,
and let $a$ be a side of $P$. We denote by $s_a$ the geodesic reflection
in $M$ about $a$, and set $P_a=P\cup s_a(P)$.

\begin{defin}  \label{admis_def}
Let $P\subset M$ be a spherical polygon. Then $P$ is {\em admissible}
if for any side $a$ the polygon $P_a$ is contained in an open
hemisphere.
\end{defin}

Note that the admissibility is a restriction only in the elliptic
case.\footnote{
We think that all of our results remain valid without this restriction.}
Unless we state otherwise, we will consider only admissible polygons,
suppressing the qualifier.

The  space $X=X(P)$ consists of
directed geodesic segments inside $P$ both of whose endpoints belong to  
$\bo P$.
Since we are using projective models, $P\subset\RR$ is a euclidean  
polygon
in any of the three cases at hand. Geodesic segments in $P$ are  
straight;
we refer to their endpoints as the {\em beginning}
$b(x)$ and the {\em end} $e(x)$. The mapping $x\mapsto(b(x),e(x))$
is an embedding $X\subset\partial P \times \partial P$;
it induces a topology on $X$.

Let $L$ be the space of oriented straight lines (rays) in the euclidean  
plane.
Endowed with the natural topology, it is a cylinder.
Besides, $L$ is a geodesic surface: geodesics in
$L$ are the pencils of rays passing through a point or the pencils
of parallel rays. Equivalently, we define this structure via the  
canonical
embedding $\RR\subset{\bf RP^2}$. Lines in $\RR$ become
points of the dual projective plane $({\bf RP^2})^*\simeq{\bf RP^2}$.  
The
geodesic surface structure of $L$ is thus induced by that of the real  
projective
plane.

Denote by $l(x)\in L$ the ray containing the chord $x$.
The mapping $X\to L$ given by $x\mapsto l(x)$ is finite-to-one;
if $P$ is convex then it is one-to-one. It induces the
structure of a geodesic surface on $X$.

Let $\Delta_+\subset X$ (resp. $\Delta_-\subset X$) be given by the
condition that $e(x)$ (resp. $b(x)$) is a corner of $P$.
Let $\Delta_0\subset X$ be the set of chords that
contain a corner of $P$ in the interior. Set
%
$$\Delta=\Delta_+\cup\Delta_-\cup\Delta_0.$$
%
By definition, $\Delta$ is a geodesic graph in $X$.
The following lemma is crucial.

\begin{lem} \label{convex_lem}
The faces of the graph $\Delta$ are geodesically convex.
\end{lem}
{\bf Proof}. Let
$x_0=A_0 B_0,x_1=A_1 B_1\in X\setminus\Delta$, and let  
$x_t=A_tB_t,\,t\in
[0,1],$
be a path in $X\setminus\Delta$ connecting them. Then the points
$A_t,\,t\in[0,1]$
(resp. $B_t,\,t\in[0,1]$) belong to the interior of a side,  
$a\subset\bo P$
(resp. $b\subset\bo P$).

Up to relabeling and orientation reversal, the points in question
may form four apriori possible configurations;
the four cases are shown in figure~\ref{figConv1}.

\begin{figure}[htbp]
\begin{center}
\input{fig101.pstex_t}
\caption{Convexity of $\Delta$}
\label{figConv1}
\end{center}
\end{figure}

Cases $3,4$ contradict to our assumptions. Indeed, in case $4$ the sides
$a,b$ intersect at their interior points. In case $3$ let $V_0,V_1$
be the endpoints of the side $b$. The two frames $A_0B_0,V_0V_1$ and
$A_1B_1,V_0V_1$ have opposite orientations, thus for some $0<t<1$
the vectors $A_tB_t,V_0V_1$ are colinear, in contradiction to the  
assumption
that
$A_tB_t\in X\setminus\Delta$.

In cases $1,2$ let $O\in\RR$ be the intersection point of the straight  
lines
$A_0B_0,A_1B_1$.\footnote{
In case $1$ the two lines may be parallel, but the same argument works.}
Let $\gamma$ be the curve in $X$ obtained by rotating the oriented line  
$A_0B_0$
about $O$ towards $A_1 B_1$ by an angle less than $\pi$.
It is a geodesic.
To show that $\gamma$ does not intersect $\Delta$,
it suffices to prove that the quadrilateral
$A_0A_1B_1B_0$ contains no vertices of $P$. Assume that it does, and let
$V$ be such a vertex. The frames $A_0V,A_0B_0$ and $A_1V,A_1B_1$
have opposite orientations. Hence, for some $0<\tau<1$ the vectors
$A_{\tau}V,A_{\tau}B_{\tau}$ are colinear, and thus the segment
$x_{\tau}=A_{\tau}B_{\tau}$
contains a corner of $P$, in contradiction to the assumptions. \qed

\medskip

We denote by $T_b:X\to X$ the billiard map and by $T_b^{-1}:X\to X$
its inverse. They are not defined on all of $X$, and there are
choices in defining them on $\Delta$. To avoid the ambiguities,
we will now precisely define  $T_b,T_b^{-1}$ on their respective  
domains of
definition.

\begin{figure}[htbp]
\begin{center}
\input{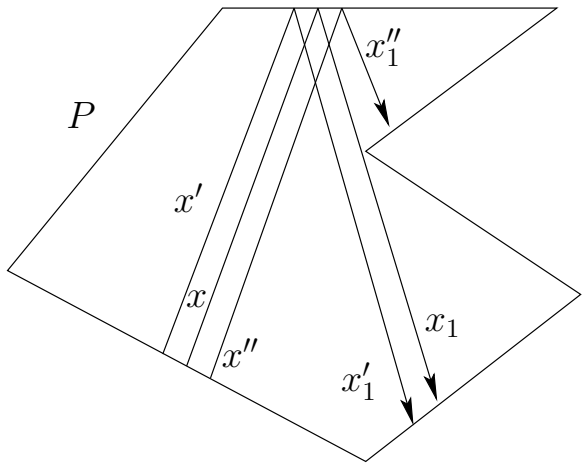}
\caption{Discontinuity of the billiard map at a point in  
$T_b^{-1}(\Delta_0)$}
\label{figDiscon}
\end{center}
\end{figure}

Let $x=[b(x),e(x)]$ be a phase point. Then $x_1=T_b(x)$ (resp.
$x_{-1}=T_b^{-1}(x)$) is not
defined iff $e(x)$ (resp. $b(x)$) is a corner of $P$.  Thus,
the natural domain of definition for $T_b$ (resp. $T_b^{-1}$) is
$X\setminus\Delta_+$ (resp. $X\setminus\Delta_-$).
Suppose now that $e(x)$ is not a corner, and let $a$ be the unique side  
of
$P$ that contains $e(x)$. Let $l_1(x)$ be the reflection
of $l(x)$ about $a$.  Let
$x_1=T_b(x)=[e(x),e(x_1)]\subset l_1(x)$ be the longest chord such that
$[e(x),e(x_1)]\in X$. This defines $T_b$ on $X\setminus\Delta_+$.
The definition of $T_b^{-1}$ on $X\setminus\Delta_-$ is analogous,
and we leave it to the reader. Set
\begin{equation}
\label{graph_eq1}
\Gamma=\Delta\cup T_b^{-1}(\Delta\setminus\Delta_-),\
\Theta=(\Delta\setminus\Delta_-)\cup T_b^{-1}(\Delta_0).
\end{equation}
\begin{prop}  \label{discont_prop}
The following holds:\\
1. The subset $\Theta\subset X$ is the set of discontinuities of
the billiard map.\\
2. The faces of $\Gamma$ are geodesically convex.
\end{prop}
{\bf Proof}. It is elementary to see that  $\Delta_+$ belongs to the  
set of
discontinuities. (The billiard map is not even defined on $\Delta_+$.)
Let $x\in\Delta_0$.  Then
arbitrarily close to $x$ there are phase points $x',x''$
such that $e(x'),e(x'')$ belong to distinct, and non-adjacent,  sides  
of $P$.
Hence $b(T_b(x'))$ and $b(T_b(x''))$ are not close to each other, and
$T_b$ is discontinuous at $x$. Let $x\notin\Delta_+$.
Then $x_1=T_b(x)$ is defined. If however $x_1\in\Delta_0$,
then arguing as above, we conclude that $e(T_b(x'))$ and $e(T_b(x''))$  
are not
close to each other, and $T_b$ is discontinuous at $x$. See
Figure~\ref{figDiscon}. On the other hand, if $x_1\notin\Delta_0$, then  
$x$
is a point of continuity of $T_b$. This proves the first claim.

Let $x',x''$ belong to a face of $\Gamma$.
Arguing as in the proof of Lemma~\ref{convex_lem}, we see that
$e(x'),e(x'')$ belong to the same side, $a$, of $P$.
Set $Q=P_a$, and let $Y=X(Q)$. The ``reflection" trick associates
with any phase point $x\in X$ such that $e(x)\in a$ a phase point $y\in Y$
\cite{Ta95}. The correspondence $x\mapsto y$ maps a curve
$x(t),\,0\le t\le 1,$ connecting  $x',x''$ in $X$ into
a curve $y(t),\,0\le t\le 1,$ connecting  $y',y''$ in $Y$.
Let $\Delta(Q)\subset Y$ be the corresponding set of phase points
containing corners.
The condition  $x(t)\in X\setminus\Gamma,\,0\le t\le 1,$
implies that $y(t)\in Y\setminus\Delta(Q),\,0\le t\le 1.$
By Lemma~\ref{convex_lem}, phase points $y',y''$ are connected by a  
geodesic,
$\beta$, in $Y\setminus\Delta(Q)$. Folding $\beta$ back into $X$ by the
geodesic reflection $s_a$, we obtain the geodesic $\alpha$
connecting $x',x''$ in $X\setminus\Gamma$.
This establishes the second claim.\qed

\vspace{3mm}

We will now cast the billiard mapping as a piecewise convex  
transformation.

\begin{thm}  \label{bil_piece_prop}
Let  $M$ be a simply connected surface of constant curvature, and
let $P\subset M$ be a geodesic polygon. Let
$X=X(P)$ be the phase space of the billiard map
endowed with the structure of a geodesic surface.
Let $T_b:X\to X$ be the billiard map.

Let $\Gamma$ be the graph defined by  (\ref{graph_eq1}),
and let $\PP:X=X_{\alpha}\cup X_{\beta}\cup \cdots$ be the corresponding
convex geodesic partition.

Then $\Gamma$ is a convex geodesic graph, and $T_b$ is continuous
on the faces of $\Gamma$. Let $T$ be the map that coincides
with $T_b$ on the open faces of $\Gamma$ and extends to
their closures by continuity. Then  $(X,T,\PP)$ is
a piecewise convex transformation; the billiard map $T_b$
and the transformation $(X,T,\PP)$ coincide on $X\setminus\Gamma$.
\end{thm}
{\bf Proof}. Most of the statements have been established
in the course of the preceding discussion. For instance,
Proposition~\ref{discont_prop}
asserts that $\Gamma$ is a convex geodesic graph. It remains to show  
that
the restriction of $T_b$ to any open face of $\Gamma$ is a geodesically  
convex
transformation.

The geodesic structure on $X$ has been defined via ray  
focusing.\footnote{
We use the self-explanatory language of geometric optics
in discussing billiard dynamics \cite{Ta95}.} Geodesic reflections send  
focusing
beams into focusing beams. Let $Y\subset X$ be a face of $\Gamma$. As we
have seen in the proof of Proposition~\ref{discont_prop}, the endpoints  
$b(x)$
of all phase points $x\in Y$ belong to the same side, $a$, of $P$.
Hence, for all $x\in Y$ the new phase point $x_1=T_b(x)$ is obtained via
the geodesic reflection $s_a$. Hence $T_b$ sends geodesic segments of  
$Y$
into geodesic segments.  \qed

\vspace{3mm}

It has been traditional to study the complexity of billiard
orbits in a polygon $P$ using the coding by sides of $P$ \cite{Ta95}.
Let $(X,T,\PP)$ be as in Theorem~\ref{bil_piece_prop}.
The complexity of billiard orbits with respect to the coding by
atoms of $\PP$ coincides with the complexity of $(X,T,\PP)$.
We will study it, and then apply our results to the traditional
billiard complexity.

\vspace{3mm}

First, we will establish a dictionary between the
language of billiard orbits \cite{Gu03,MT02,Ta95} and that
of piecewise convex transformations.
A billiard orbit of (combinatorial) length $n$ is a sequence
$\gamma= x_1,\dots,x_n$ of phase points  such that $x_{i+1}=T_b(x_i)$.
Geometrically, $\gamma$ is a sequence of $n$ consecutive chords of $P$,
where each chord is the reflection of the preceding one.
In particular, the points $e(x_1),\dots,e(x_{n-1})$ are not corners.
We will also say that $\gamma$ is a {\em $n$-segment orbit}.
If $e(x_n)$ is not a corner, then the {\em $1$-step forward extention}
is $x_1,\dots,x_n,x_{n+1}$, where $x_{n+1}=T_b(x_n)$.
If $b(x_1)$ is not a corner, then the {\em $1$-step backward extention}
is $x_{-1},x_1,\dots,x_n$, where $x_{-1}=T_b^{-1}(x_1)$.
We can iterate these extentions in obvious ways.

A billiard orbit is {\em regular} if it does not
contain corners. Otherwise, $\gamma$ is {\em singular}.
If $\gamma$ does not contain corners of $P$ in its interior,
but the endpoints $b(\gamma)=b(x_1),e(\gamma)=e(x_n)$ are corners,
we say, following \cite{Ka87}, that $\gamma$ is a {\em generalized  
diagonal}
of length $n$.

The notion of generalized diagonals  works well for the
billiard in a convex polygon. We introduce an extension of this notion
which works for arbitrary polygons.
\begin{defin}  \label{sing_def}
An $n$-segment billiard orbit $\gamma= x_1,\dots,x_n$ is {\em strongly  
singular}
if $x_2,\dots,x_{n-1}$ do not contain corners but $x_1$ and $x_n$ do.
\end{defin}

A {\em family of billiard orbits} is a one-parameter family
$\gamma(t)=x_1(t),\dots,x_n(t)$ where the mapping $t\mapsto\gamma(t)$
is injective and continuously differentiable.
In particular, we will consider families of strongly singular billiard  
orbits.
Any such is contained in a unique maximal family, and
we will consider only them, suppressing the qualifier {\em maximal}.
A strongly singular billiard orbit is {\em isolated} if
it is not contained in a family of such orbits.

\begin{prop}   \label{transition_prop}
Let $P\subset M$ be a geodesic polygon, and let $(X,T,\PP)$ be the
corresponding piecewise convex transformation. Let $n\ge1$. Then:

\noindent 1. There is a bijection between the set $V_{\ess}(n)$ and the  
union of
the
sets of isolated $(n+2)$-segment and $(n+3)$-segment strongly
singular billiard orbits;

\noindent 2. There is a bijection between the set $OE(n)$ and the union  
of the
sets of families of $(n+2)$-segment  and $(n+3)$-segment strongly
singular billiard orbits.
\end{prop}
{\bf Proof}. Let $x\in X\setminus\Gamma_n$. Iterating the billiard map,
we obtain an $(n+1)$-segment billiard orbit $x=x_1,\cdots,x_n,x_{n+1}$,
where $x_{n+1}=T_b^n(x)$. Let $x_{n+2}=T_b(x_{n+1})$. Then  
$x\in\Gamma_{n+1}$
iff either $x_{n+1}\in\Delta_0\cup\Delta_+$ or $x_{n+1}\in  
X\setminus\Delta$
but $T_b(x_{n+1})=x_{n+2}\in\Delta$.

Set $x_{-1}=T_b^{-1}(x)$. Then $x\in\Gamma_{-1}$ iff
$x_{-1}\in\Delta\setminus\Delta_+$.
Thus, we have obtained a surjective map from
$<\Gamma_{-1}\vee(\Gamma_{n+1}\setminus\Gamma_n)>\subset X$ to the  
union of
set of $(n+2)$-segment and $(n+3)$-segment strongly singular billiard  
orbits.
(The orbit corresponding to $x$ is $x_{-1},x_1,\cdots,x_n,x_{n+1}$
in the former case, and
$x_{-1},x_1,\cdots,x_n,x_{n+1},x_{n+2}$  in the latter.)

The map above is a bijection.
Moreover, isolated orbits correspond to the  essential vertices of
$\Gamma_{-1}\vee(\Gamma_{n+1}\setminus\Gamma_n)$, and families of orbits
correspond to
overlapping edges of $\Gamma_{-1}$ and  
$\Gamma_{n+1}\setminus\Gamma_n$.\qed

\section{Complexity of the billiard in a polygon}  \label{cases}
We will now apply the preceding material to the complexity of  billiards
in geodesic polygons on surfaces of constant curvature. In the  
beginning of
this section we consider  the three cases at hand simultaneously,
emphasizing their similarities.
\subsection{Arbitrary curvature, any polygon}
\label{arbit_sub}
\begin{figure}[htbp]
\begin{center}
\input{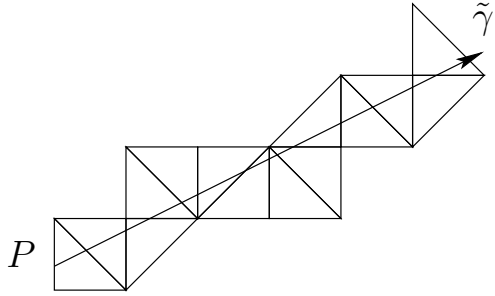}
\caption{Unfolding  a billiard trajectory}
\label{fig105}
\end{center}
\end{figure}

Let $M$ be a surface of constant curvature $\chi=0,\pm 1$, and let
$P\subset M$ be a geodesic polygon. The operation of {\em unfolding}
sends billiard orbits in $P$ into geodesics in $M$ \cite{Ta95}, see figure \ref{fig105}.
If $\gamma=x_1,x_2,\dots,x_n$ is a billiard orbit,
then its unfolding is the geodesic
$\tilde{\gamma}=x_1,\tilde{x}_2,\dots,\tilde{x}_n$
where the segments $x_i,\tilde{x}_i$ differ by an isometry of $M$.

Let $G(P)\subset Iso(M)$ be the group generated by
the geodesic reflections in the sides of $P$. Denote by  
$G^{(n)}(P)\subset G(P)$
the set of elements obtained from products of at most $n$ reflections.
(These are the elements of length at most $n$.)
Then $G^{(n)}(P),0\le n < \infty,$ is an increasing tower of finite  
sets;
their union is $G(P)$.

\begin{lem}  \label{conj_lem}
1. Let $\chi\le 0$. Then all strongly singular billiard orbits in $P$  
are
isolated.

\noindent 2. Let $\chi> 0$. Let $\gamma=x_1,x_2,\dots,x_n$ be a
strongly singular billiard orbit in $P$ and let
$\tilde{\gamma}=x_1,\tilde{x}_2,\dots,\tilde{x}_n$ be its unfolding.
Suppose that $\gamma$ extends to a family of strongly singular billiard  
orbits.
Then there are corner points $p_1\in x_1,p_n\in x_n$ and an isometry  
$g$ of $M$
such that
$g(x_n)=\tilde{x}_n$ and the points $p_1,g(p_n)$ either coincide or
are antipodal.
\end{lem}
{\bf Proof}. Let
$\gamma(t)=x_1(t),x_2(t),\dots,x_n(t),\,0\le t \le 1,$ be a family of
strongly singular billiard orbits. The segments $x_1(t)$ (resp.  
$x_n(t)$)
pass through the same corner point $p_1$ (resp. $p_n$). We view
$\gamma(t)$ as a beam of trajectories that emanate from the focusing
point $p_1$ and refocus at $p_n$.

The unfolding $\gamma\mapsto\tilde{\gamma}$ transforms directed  
billiard orbits
into directed geodesics in $M$, preserving the length and sending
focusing beams of billiard orbits
into focusing beams of geodesics\cite{Ta95}.\footnote{
Although \cite{Ta95} assumes that
$\chi=0$, the argument applies to arbitrary $\chi$.}
The beam $\tilde{\gamma}(t)$
emanates from $p_1\in M$ and refocuses at $g(p_n)\in M$
where $g\in G(P)$. The distance between $p_1,g(p_n)$ along any
geodesic $\tilde{\gamma}(t)$ is the same as the distance
between $p_1,p_n$ along the orbit $\gamma(t)$, hence positive.

If  $\chi\le 0$, this implies $p_1\ne g(p_n)$; thus
the geodesic beam $\tilde{\gamma}(t)$ has two different focusing points.
This is impossible, which proves the first claim.
If $\chi> 0$, i.e., $M$ is the sphere, the focusing
points in question are either antipodal or they coincide.\qed

\vspace{3mm}

We denote by $s(k)$ (resp. $fs(k)$) the number of isolated (resp.  
families of)
strongly singular $k$-segment billiard orbits in $P$.
Set
\begin{equation}  \label{strong_sing_eq1}
S(n)=\sum_{3 \le k \le n}s(k),\,FS(n)=\sum_{3 \le k \le n}sf(k).
\end{equation}
\begin{prop}     \label{str_sing_prop}
Let $P$ be a geodesic polygon on a surface $M$ of constant curvature  
$\chi$.
Let $(X,T,\PP)$ be the associated piecewise convex transformation,
and let $f(\cdot)$ be the billiard complexity corresponding to the  
partition
$\PP$.
Then there are positive integers $q_1,q_2$ depending on $P$, so that:

\noindent 1. If $\chi=1$, then
\begin{equation}  \label{compl_polyg_eq}
f(n) = q_1 + q_2n + S(n+1) -  FS(n+1)    +
2\left(\sum_{3\le k \le n}S(k) - \sum_{3\le k \le n}FS(k)\right).
\end{equation}
\noindent 2. If $\chi\le 0$, then
\begin{equation}  \label{compl_curv_eq}
f(n) = q_1 + q_2n + S(n+1) + 2\sum_{3\le k \le n}S(k).
\end{equation}
\end{prop}
{\bf Proof.} We use Proposition~\ref{transition_prop} to express
the right hand side of (\ref{complex_eq}) via
the numbers of (families of) strongly singular billiard orbits.
This yields the first claim.
Taking into account Lemma~\ref{conj_lem}, we obtain the second claim.
The factor of $2$
in (\ref{compl_polyg_eq},~\ref{compl_curv_eq}) is due to
the fact that the contribution of $k$-segment strongly singular orbits
with $k\le n$ is counted twice. See Proposition~\ref{transition_prop}.
  \qed
\vspace{3mm}

There is a correspondence between strongly singular
billiard orbits and generalized diagonals. Let
$\gamma=x_1,x_2,\dots,x_{n-1},x_n,$
$n\ge 3,$
be a strongly singular orbit. Let $p\in x_1$ (resp. $q\in x_n$) be the  
last
(resp. the first) corner point. Note that $p\ne e(x_1), q\ne b(x_n)$.
Set $\tilde{x}_1=[p,e(x_1)],\tilde{x}_n=[b(x_n),q]$.
Then $\tilde{\gamma}=\tilde{x}_1,x_2,\dots,x_{n-1},\tilde{x}_n$
is a generalized diagonal. The mapping $\gamma\mapsto\tilde{\gamma}$
is a bijection, and it preserves the combinatorial length.\footnote{
However, it does not preserve the property to be isolated.}

Let $gd(n)$ (resp. $fgd(n)$) be the number of $n$-segment isolated
(resp. families of) generalized diagonals, $n\ge 3$. Set
\begin{equation} \label{defin_eq}
GD(n)=\sum_{3 \le k \le n}gd(k),\ FGD(n)=\sum_{3 \le k \le n}fgd(k).
\end{equation}
\begin{thm}     \label{gen_diag_thm}
Let $P$ be a geodesic polygon on a surface $M$ of constant nonpositive
curvature. Let $(X,T,\PP)$ be the associated piecewise convex  
transformation,
and let $f(\cdot)$ be the billiard complexity corresponding to the  
partition
$\PP$.
Then there are integers $c_1,c_2$ depending on $P$, so that:
\begin{equation}  \label{asymp_eq1}
f(n) = c_1 + c_2n + GD(n+1) + 2\sum_{3\le k \le n}GD(k).
\end{equation}
\end{thm}
{\bf Proof}. By Lemma~\ref{conj_lem}, strongly singular billiard orbits  
are
isolated. The same argument shows that generalized diagonals are  
isolated,
as well. The claim now follows from (\ref{compl_curv_eq}) and the
correspondence
between strongly singular billiard orbits and generalized diagonals.   
\qed

\vspace{3mm}

We will now apply the preceding material to the complexity, $F(\cdot)$,
of the coding billiard orbits by the sides of $P$. To avoid  
ambiguities, we
will define this coding now. Let $P$ have $p$ sides $a,b,\dots,$ and
denote by $X(a,b)\subset X$ the set of phase points $x$ such that
$b(x)\in a,e(x)\in b$. As $a\ne b$ run through the pairs of sides
such that $X(a,b)$ has a nonempty interior, the sets $X(a,b)$ form the  
partition
$\PP_{\mbox{side}}$ of the phase space. Then the {\em side complexity},
$F(\cdot)$, is the complexity of the billiard map with respect to
$\PP_{\mbox{side}}$.

It will be useful to have a direct definition of the side complexity.
We denote by $\LL_{\mbox{side}}$ the full language on the alphabet of  
sides
of $P$. Let now $\gamma=x_1,\dots,x_n$ be a regular $n$-segment  
billiard orbit.
Then $b(x_1),\dots,b(x_n),e(x_n)$ are interior points of sides
$s_1,\dots,s_n,s_{n+1}$, and we set
$$\sigma_{\mbox{side}}(\gamma)=s_1,\dots,s_n,s_{n+1}\in\LL_{\mbox{side}} 
(n+1).$$
The language $\Sigma_{\mbox{side}}$ is the range of the mapping, and we  
set
$\Sigma_{\mbox{side}}(n)=\Sigma_{\mbox{side}}\cap\LL_{\mbox{side}}(n)$.
Then the {\em side complexity} satisfies
\begin{equation}  \label{side_eq1}
F(n) = |\Sigma_{\mbox{side}}(n+1)|.
\end{equation}

To compare complexities, we will use an elementary  lemma.
\begin{lem}  \label{compar_lem}
Let $(X,T,\PP)$ be a piecewise convex transformation, and let $f(\cdot)$
be the corresponding complexity. Let $\PP_1$ be a partition, such that
$\PP\prec\PP_1$, and let $f_1(\cdot)$ be the complexity of $(X,T,\PP)$  
with
respect
to $\PP_1$. Then $f_1(n)\le f(n)$.
\end{lem}
\begin{corol} \label{side_cor}
Let $P$ be a geodesic polygon on a surface $M$ of constant
curvature $\chi$, and let $F(\cdot)$ be the complexity
of coding billiard orbits by the sides of $P$.
Then  the following inequalities hold:

\noindent 1. If $\chi=1$, then there are constants $q_1,q_2$ such that
\begin{equation}  \label{side_eq2}
F(n) \le q_1 + q_2n + S(n+1) -  FS(n+1) +
2\left(\sum_{3\le k \le n}S(k)-\sum_{3\le k \le n}FS(k)\right).
\end{equation}

\noindent 2. If $\chi\le 0$, then there are constants $c_1,c_2$ such that
\begin{equation}  \label{side_eq3}
F(n) \le  c_1 + c_2n + GD(n+1) + 2\sum_{3\le k \le n}GD(k).
\end{equation}
\end{corol}
{\bf Proof}. Let $\PP$ be the defining partition of the
associated piecewise convex transformation.
By Proposition~\ref{discont_prop}, the partitions  
$\PP,\PP_{\mbox{side}}$
satisfy $\PP\prec\PP_{\mbox{side}}$.
Now the first (resp. second) claim is immediate from
Proposition~\ref{str_sing_prop}
(resp. Theorem~\ref{gen_diag_thm})
and Lemma~\ref{compar_lem}.  \qed

\subsection{Arbitrary curvature, convex polygon}
\label{convex_sub}
The preceding considerations drastically simplify for
convex polygonal billiard tables.
\begin{prop}  \label{convex_prop}
Let $P\subset M$ be a convex geodesic polygon on a surface
of constant curvature $\chi$. Let $F(\cdot)$ be
the traditional complexity of billiard orbits in $P$
coming from the coding of billiard orbits in $P$ by the sides they hit.
Then there exist constants $c_1,c_2$ depending only on
$P$ such that the following holds:

\noindent 1. Let $\chi=1$. Then
\begin{equation}  \label{convex_eq1}
F(n) =  c_1 + c_2n + \sum_{3\le k \le n}GD(k) - \sum_{3\le k \le  
n}FGD(k);
\end{equation}

\noindent 2. Let $\chi\le 0$. Then
\begin{equation}  \label{convex_eq2}
F(n) =  c_1 + c_2n + \sum_{3\le k \le n}GD(k).
\end{equation}
\end{prop}
{\bf Proof}.
In the notation of Section~\ref{innerbill}, we have
$$
\Delta_0=\emptyset,\ \Theta=\Delta_+,\  
\Gamma=\Delta=\Delta_+\cup\Delta_-.
$$
The standard defining partition is $\PP=\PP(\Gamma)$, and  
we have $\PP_{\mbox{sides}}=\PP(\Delta)$.
In view of the above, $\PP_{\mbox{sides}}=\PP$, and hence the two
complexities coincide: $F(n)=f(n)$.

For a convex polygon, the notion of strongly singular
billiard orbits and the notion of generalized diagonals coincide.
The argument of Proposition~\ref{transition_prop} works. Due to
 the coincidences we just pointed out, it
establishes a bijection between $V_{\mbox{ess}}(n)$ (resp. $OE(n)$)
and the  set of isolated (resp. families of)
$n$-segment generalized diagonals in the convex polygon $P$.

From this bijection and (\ref{complex_eq}), we obtain the first  
 claim. Now the first statement of  Lemma~\ref{conj_lem} yields the second  
claim.  \qed

\begin{rem} \label{low-order_rem}
{\rm The preceding proof yields expressions for the constants
$c_1,c_2$. From (\ref{complex_eq}), we have
$c_1=2|\PP_1|-|\PP_2|,c_2=|\PP_2|.$ Let the polygon $P$ have $p$ sides.
The integer $|\PP_1|$ (resp. $|\PP_2|$) is the number
of types of $1$-segment (resp. $2$-segment) billiard orbits
in $P$. By convexity, $|\PP_1|=p(p-1)$. However, $|\PP_2|$
is not determined by $p$ alone.}
\end{rem}

We will now apply the preceding material to the complexity
of polygonal billiards. We will consider the three cases
separately.

\subsection{The euclidean case}  \label{euclidin_sub}
There is a considerable literature on the billiard dynamics in euclidean
polygons. See \cite{Ta95} for references. Many basic questions remain  
open
\cite{Gu03}. One of them is whether the complexity (of the coding by  
sides)
of billiard orbits grows (at most) polynomially.
It is known that the growth is subexponential  
\cite{Ka87,GKT,GH97},\footnote{
The discussion in \cite{Ka87,GKT} is restricted to simply connected polygons.}
implying that the billiard in a
euclidean polygon has zero topological entropy. By  
Corollary~\ref{side_cor},
the complexity is bounded from above by the counting function for  
generalized
diagonals. The latter also grows subexponentially \cite{GH97},
and is believed to grow at most polynomially \cite{Gu03}.

A euclidean  polygon is {\em rational} if all of its angles are rational
multiples of $\pi$. Rational polygons play an important role in the  
subject
\cite{MT02,Ta95}.
\begin{thm}  \label{rat_complex_thm}
Let $P$ be a rational euclidean polygon, and let $F(\cdot)$
be the side complexity of billiard orbits in $P$. Then there exists  
$c>0$
such that
\begin{equation}  \label{cub_eq0}
F(n) < cn^3 .
\end{equation}
If $P$ is convex, then there exist positive constants
$c_1,c_2$ such that
\begin{equation}  \label{cub_eq1}
c_1n^3 < F(n) < c_2n^3 .
\end{equation}
\end{thm}
{\bf Proof.} By a theorem of H. Masur \cite{MT02},  $GD(n)$ grows  
quadratically
for rational polygons. More precisely, there exist
positive constants $c_1',c_2'$ such that
\begin{equation}  \label{gen_diag_ess_vert}
c_1'n^2 < GD(n) < c_2'n^2.
\end{equation}
The first claim follows from this and the bound (\ref{side_eq3}).
The second claim follows the same way from the formula above
and the identity (\ref{convex_eq2}).\qed

\begin{figure}[htbp]
\begin{center}
\input{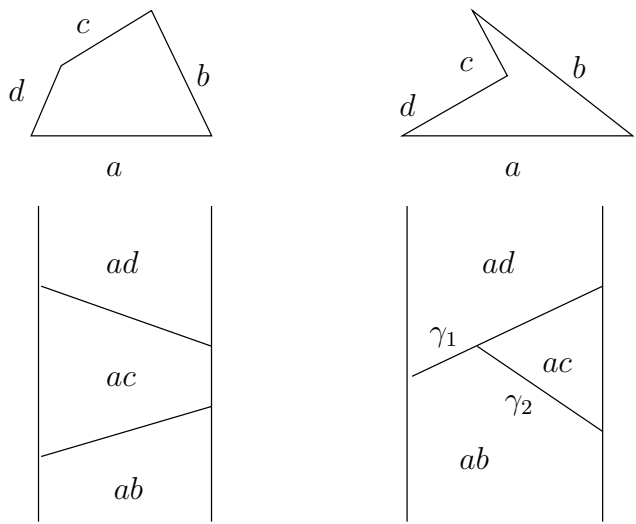}
\caption{Singular set  for the billiard in a euclidean quadrilateral:
i) convex; ii) nonconvex}
\label{fig2}
\end{center}
\end{figure}

\begin{rem} \label{convex_rem}
{\rm We believe that the cubic lower bound on complexity  
(\ref{cub_eq1})
is valid for arbitrary rational polygons. However, the geometry
of the billiard map is much simpler in the
convex case, as we saw in section~\ref{convex_sub}.
Figure~\ref{fig2} shows that the singular graph of
the billiard map in a quadrilateral is much more
complicated in the nonconvex case. The bounds (\ref{cub_eq1})
were obtained in \cite{CHT03}.}
\end{rem}

\vspace{3mm}

\subsection{The elliptic and the hyperbolic geometries}     
\label{spherein_sub}
The billiard in a spherical polygon has two salient features.
First, the billiard map is a piecewise isometry. Second,
a spherical polygon may have non-isolated strongly singular orbits,
in particular, non-isolated generalised diagonals.

We will use the notation $Geo\,(M)$ for the space of oriented geodesics
on $M$. Let $\phi:Geo\,(S^2) \to S^2$ be the standard diffeomorphism.
(See Section~\ref{neg_curv_sub} for a discussion
of spherical duality.) Pulling back by  $\phi$ the round metric on  
$S^2$,
we obtain an invariant metric; the distance between two great
circles is the angle between them. Metric geodesics and geodesics
of the projective structure coincide.

Let now $P\subset S^2$ be a polygon.
Let $X=X(P)$ be the billiard map phase space, and let $\varphi:X \to  
Geo\,(S^2)$
be the map introduced in Section~\ref{innerbill}.
Pulling back by  $\varphi$ induces not only the structure of a piecewise
convex surface on $X$ but also a metric on convex pieces of the phase  
space.
Locally, the billiard map $T_b:X\to X$ is a geodesic reflection, i.e.,  
an
isometry.
We summarize this discussion as a proposition.
\begin{prop}  \label{piec_isom_prop}
Let $P$ be an arbitrary spherical polygon. Then the billiard map phase  
space
is a piecewise metrically convex Riemannian surface. The billiard map
is a piecewise isometry.
\end{prop}

This observation has consequences for complexity of billiard orbits.
\begin{thm}  \label{spheresubexp}
Let $P$ be any spherical polygon.\footnote{
We do not assume that $P$ is an admissible polygon.}
Then the complexity of the coding
of billiard orbits by sides of $P$ grows  subexponentially.
\end{thm}
{\bf Proof.} Let $(X,T,\PP)$ be the associated piecewise convex  
transformation.
Denote by $f(\cdot)$ its complexity, and let $F(\cdot)$ be the
standard complexity of billiard orbits. By the proof of
Corollary~\ref{side_cor},
$F(n)\le f(n)$. Hence, it suffices to show that $f(n)$ grows  
subexponentially.
By Proposition~\ref{piec_isom_prop},
$(X,T,\PP)$ is a piecewise isometry on a convex partition.
By Theorem 4.2 of \cite{GH97},
its complexity has subexponential growth.   \qed

\vspace{3mm}

The examples below illustrate peculiarities of spherical polygonal
billiards.

\medskip

\noindent{\bf Example 1}. Let $P\subset S^2$ be a polygon such that $G(P)$ is
a finite group. For instance, $P$ may be the fundamental domain
of a finite, generated by reflections group of isometries.\footnote{
These polygons are well known. See, e.g., \cite{GN04}.}
Then every billiard orbit  in $P$ is periodic. Moreover, there is
a finite number of symbolic codes corresponding to prime periodic
orbits, hence complexity is bounded.

\medskip

\noindent{\bf Example 2}. Let $P$ be a ``bigon"; it is bounded by two
geodesics, $a,b$ connecting the North and the South poles.
Let $\alpha$ be the angle  between them.\footnote{
Note that $P$ is not an admissible polygon.}
If $\alpha$ is $\pi$-rational, then we are in the situation
of Example 1. We will now discuss the case when
$\alpha$ is $\pi$-irrational.

First, we point out that in any case the set $\Sigma_{\mbox{side}}(n)$
consists of $2$ elements: $a,b,a,b,\dots$ and $b,a,b,a,\dots$.
Thus, $F(n)=2$.

There is an obvious periodic orbit.
It corresponds to the intersection of $P$ with the equator, and it has
$2$ segments, perpendicular to $a,b$. We will denote this orbit
by $\gamma_0$.

\noindent{\bf Claim.} Let $\alpha$ be $\pi$-irrational.
Then $\gamma_0$ is the only prime periodic orbit in $P$.

We will show that any periodic orbit $\gamma$ is a multiple of
$\gamma_0$. We can assume that $\gamma$ has an even number, $2m$,
segments, and that its symbolic code is $b,a,\dots,b,a$. Then
the element, $g(\gamma)$, of the group $G(P)$ that we obtain by tracing  
$\gamma$
is $(\rho)^m$ where $\rho$ is the rotation about the vertical axis
by the angle $2\alpha$.

Let $\ell(\gamma)$ be the spherical geodesic corresponding to
$\gamma$. (Note that $\ell(\gamma)$ differs, in general, from the
unfolding $\tilde\gamma$, which is a geodesic segment along
$\ell(\gamma)$.) The periodicity of $\gamma$ implies that
$\ell(\gamma)$ is invariant under $g(\gamma)$, which rotates the sphere  
about
the vertical axis by $2m\alpha$. The only geodesic invariant
under this (nontrivial!) rotation is the equator, which implies the  
claim.

\medskip

By convention, a periodic billiard orbit in $P$ does not pass
through its corners. In particular, it cannot trace the boundary
of $P$. It is not known if every euclidean polygon has a periodic orbit
\cite{Gu03}. Below we present a spherical polygon without periodic  
orbits.

\medskip

\noindent{\bf Example 3}. For $0 < \alpha < 2\pi$ let $Q=Q(\alpha)$ be  
the
isosceles spherical triangle with two right angles, and whose third
angle is $\alpha$. If $\alpha$ is $\pi$-rational then every billiard  
orbit
in $Q$ is periodic. If $\alpha$ is $\pi$-irrational, then $Q$ has no
periodic billiard orbits. We outline a proof below.

The triangle $Q$ is obtained from the bigon $P$ of Example 2 by folding
it about the equator. In this situation, every billiard orbit,
$\gamma$, in $Q$ uniquely lifts to a billiard orbit $\tilde\gamma$ in  
$P$;
the orbit $\gamma$ is periodic iff so is $\tilde\gamma$. If $\alpha$ is
$\pi$-rational,
then $Q$ satisfies the conditions of Example 1.
Let $\alpha$ be $\pi$-irrational, and let $\gamma$ be a periodic
orbit in $Q$. By preceding remark and Example 2, $\tilde\gamma$
runs along the equator. Thus, $\gamma$ traces the boundary
of $Q$.

\vspace{5mm}

We will now discuss the hyperbolic case. A positive function,  
$s(\cdot)$,
of natural argument is {\em subexponential} if
$s(n)< e^{hn},s^{-1}(n)< e^{hn}$ for any $h>0$ and all sufficiently great $n$. 
\begin{thm}  \label{hyperb_exp}
Let $P\subset\HT$ be a geodesic polygon, and let $f(\cdot)$
be the complexity of the coding
of billiard orbits by sides of $P$. Denote by $h_{\rm{top}}$
the {\em topological entropy} of the billiard map in $P$.\\
Then $h_{\rm{top}}>0$; there exists a subexponential function
$s(\cdot)$  such that $f(n)=s(n)e^{h_{\rm{top}}n}$.
\end{thm}
{\bf Proof.} The billiard flow of $P$ is (uniformly) hyperbolic
\cite{GGS99}.\footnote{
This is a very special
case of a general result in \cite{GGS99}. It can be obtained
directly via standard techniques. Although specialists seem to be aware  
of it,
to our knowledge this is not in the literature.}
Thus, the metric entropy of the billiard flow with respect to the
Liouville measure is positive. By Abramov's formula, the metric entropy  
of the
billiard map in $P$ is positive, as well. The metric entropy
is a lower bound on the topological entropy, hence our first claim.

Let $\PP,\PP_{\mbox{side}}$ be our defining partition and
the partition by sides of the phase space of the billiard map.
Let $\alpha(t),\beta(t)$ be
infinite billiard orbits that
bounce of the same sides of $P$ as $-\infty < t < \infty$.
Let $\tilde{\alpha}(t),\tilde{\beta}(t)$ be their unfoldings.
Then  $\tilde{\alpha}(t),\tilde{\beta}(t)$ are infinite
geodesics in $\HT$ such that the distance between them is bounded
for  $-\infty < t < \infty$. Hence
$\tilde{\alpha}(-\infty)=\tilde{\beta}(-\infty)$
and $\tilde{\alpha}(\infty)=\tilde{\beta}(\infty)$ implying
$\tilde{\alpha}=\tilde{\beta}$. Thus, $\alpha=\beta$,
and hence $\PP_{\mbox{side}}$ is a generating partition.

Since $\PP\prec\PP_{\mbox{side}}$, the convex partition
$\PP$ is generating as well. By \cite{GH97},
the complexity of $(X,T,\PP)$ has the form
$t(n)e^{h_{\rm{top}}n}$, where $t(\cdot)$ is a subexponential
function. The relation $\PP\prec\PP_{\mbox{side}}$ implies
$f(n)\le t(n)e^{h_{\rm{top}}n}$. The cardinality of the number of
atoms of $\PP(n)$ that
partition an atom of $\PP_{\mbox{side}}(n)$ grows at most polynomially;
thus $f(n)\ge t(n)n^{-d}e^{h_{\rm{top}}n}$ for some positive integer  
$d$. \qed

\begin{rem} \label{lorenz_rem}
{\rm   Unlike $Geo\,(S^2)$, the spaces
$Geo\,(\HHH^2)$ and $Geo\,(\R^2)$
do not have Riemannian metrics, invariant under the natural actions
of the groups of
isometries of  $\HHH^2$ and $\R^2$, respectively. However, there exists
an  invariant Lorentz metric on $Geo\,(\HHH^2)$. We describe it below.

Let $H$ be the upper sheet of the
hyperboloid $z^2 -x^2 -y^2 =1$ in $\R^3$ equipped with the
pseudo-Riemannian metric
$g=d x^2 + d y^2 - d z^2$. The induced metric on $H$ is the metric of
constant negative
curvature. Isometries of $H$ are the restrictions of $g$-orthogonal
transformations of the ambient space. Geodesics in
$H$ are its intersections with planes through the origin.

Given such a plane, its
$g$-orthogonal complement is a line that intersects the hyperboloid
of one sheet
$H_1=\{z^2 -x^2 -y^2 =-1\}$ at two antipodal points; if the plane is
oriented then so is
the line, and one can canonically choose one of these intersection  
points. This
construction identifies $Geo\,(\HHH^2)$ with $H_1$. The metric
$g$ induces a pseudo-Riemannian metric of signature $(1,1)$ on $H_1$.
The identification above yields a Lorentz metric on $Geo\,(\HHH^2)$.

Let now $P\subset\HHH^2$ be a geodesic polygon, and let
$X$ be the phase space of the billiard map in $P$. The canonical
geodesic space structure in $X$ is induced by the natural mapping
$p:X\to Geo\,(\HHH^2)$. Pulling back the Lorenz metric on  
$Geo\,(\HHH^2)$,
we obtain a metric $h$ on $X$. By construction, $h$ is locally  
invariant under
the billiard map.\footnote{In general, $h$ is a Lorentz metric with
singularities. If $P$ is a convex polygon,
then $h$ is regular.} Thus, the billiard map
in a  hyperbolic polygon is a piecewise Lorentz isometry.
We do not know of any applications of this observation.}
\end{rem}

\section{Complexity of polygonal outer billiards}  \label{outerbill}
Let $(X,T,\PP)$ be a piecewise convex transformation.
If the mappings $T_i:P_i\to X$ are isometries, then
$(X,T,\PP)$ is a {\em piecewise convex isometry}.
Piecewise isometries in  one dimension are the interval exchange maps.
They arise, in particular, from the billiard in rational polygons
\cite{MT02}, and have been much studied.  We will investigate the  
complexity of
a particular class of
piecewise convex isometries -- the outer billiard transformations.

Let $M$ be a simply connected  surface of constant curvature  
$\chi=0,\pm 1$.
For $x\in M$ let $T_x:M\to M$ be the geodesic symmetry about $x$.
Let $P\subset M$ be a convex polygon with $p$ vertices $a,b,c,\dots$
listed counterclockwise.  If $\chi\ne 1$,  set $X=X(P)=M-P$.
If $\chi=1$, i. e., $M$ is the sphere, we
assume that $P$ is contained in a hemisphere.\footnote{
This will be our standing assumption; we will not restate it.}
Let $P'$ be the  antipodal polygon, and set $X=X(P)=M-P-P'$.

\begin{figure}[htbp]
\begin{center}
\input{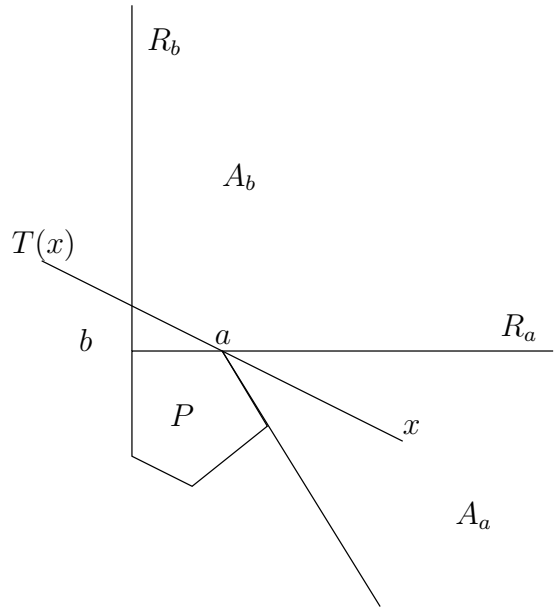}
\caption{Definition of the outer billiard map}
\label{fig3}
\end{center}
\end{figure}

For a vertex, say $a$, of $P$, let
$R_a\subset X$ be the geodesic ray  extending the side $ab$ in the
direction of $a$. The rays $R_a,R_b,\dots$ partition $X$ into convex
polygons $X_a,X_b,\dots$. See figure \ref{fig3}. We denote this
partition by  $\PP$. The statement below defines the protagonist
of this section.

\begin{defin}  \label{outer_bil_piece_prop}
Let $M$ be a simply connected  surface of constant curvature $\chi$, and
let $P\subset M$ be a convex polygon.
Set $X=M-P$ if $\chi\ne 1$, and $X=M-P-P'$ if $\chi = 1$.

Let $\PP:X=X_a\cup X_b\cup\cdots$ be the partition of
$X$ defined above. Then $\PP$ is a convex geodesic
partition.  The piecewise convex isometry $(X,T,\PP)$ defined
by the geodesic symmetries $T_a:X_a\to X,T_b:X_b\to X,\dots$
is the {\em outer billiard} about $P$, and the
geodesic surface $X=X(P)\subset M$ (it is a topological
annulus) is the {\em phase space of the outer billiard}.
\end{defin}

We will use the notation $T_o:X\to X$ for the outer billiard,
suppressing the subscript if no confusion arises.
By {\em complexity of the outer billiard} we mean the complexity
of $(X,T,\PP)$ with respect to the partition $\PP$. Now we introduce
notation and terminology that will be used throughout this section.
If $g(n),h(n)$ are two positive sequences then we write
$g \prec h$ if there is a constant $C$ such that for all $n$  
sufficiently
large $g(n) \leq C h(n)$. If $g \prec h$ and $h \prec g$, then we
write $g \sim h$; we will say that the sequences {\em have the same  
growth}
or are in the {\em same (growth) class}.
If $g \prec n^d$ then we say that $g$ grows at most polynomially
with degree $d$, or that $g$ is bounded by $n^d$.

If $G$ is a group with a finite set $S=\{s_1,\dots,s_p\}$ of generators,
we denote by $G_S^{(n)}\subset G$ the set of elements that can be  
represented
by products of at most $n$ elements of $S$ and their inverses.
The growth class of the sequence $g_S(n)=|G_S^{(n)}|$ does not depend on
the choice of $S$ \cite{dH}. If $g_S(n)\sim n^d$,
then we say that the group $G$ {\em grows polynomially, with degree  
$d$}.

Let $G=G(P)\subset Iso(M)$ be the group generated by
the geodesic reflections  in the vertices of the polygon $P$.
We will relate the growth of $G$
and the complexity of the outer billiard about $P$. We proceed to
study separately the three cases at hand.

\subsection{The euclidean case}   \label{euclid_sub}
We will obtain polynomial  bounds on the complexity of outer billiard.

\begin{thm}   \label{upperbound}
Let $P$ be a convex euclidean $p$-gon, and let $f(\cdot)$ be the   
complexity
of the outer billiard about $P$. Then $n \prec f(n)\prec n^{p+1}$.
\end{thm}
{\bf Proof.} The edges of the graph $\Gamma_n$ are parallel to the  
sides of
$P$; each edge is a segment or a half-line. Assume, for simplicity of
exposition,
that $P$ has no parallel sides. Then there are $p$ directions. For each
direction there
are $n$ parallel half-lines, hence their total number is $pn$. Since  
they
partition  $X$ into
$pn$ components,  the number of faces of $\Gamma_n$ is at least $pn$.
This yields the linear lower bound on complexity.\footnote{
We conjecture that there is a universal quadratic lower bound.}
Now, for the upper bound.

\medskip

Let $G=G(P)$, and let $S=\{T_1,\dots,T_p\}$ be its natural
set of generators. We will need a few lemmas.

\begin{lem}  \label{growthgroup}
The growth of $G$ is bounded by $n^{p-1}$.
\end{lem}
{\bf Proof.} The subgroup $H\subset G$ generated by
$T_1T_p, T_2T_p,\dots, T_{p-1}T_p$ is a quotient  group of
$\Z^{p-1}$; hence its growth is bounded by $n^{p-1}$. Since $H$
is a normal  subgroup of $G$ of index $2$, the two groups
have the same growth. \qed

Let $\Gamma=\bo\PP$, and let $\Gamma_1, \Gamma_2, \dots$ be the
canonical sequence of graphs; see Section~\ref{basic}. Let
$\gamma_n$ be the set of edges of $\Gamma_n\setminus\Gamma_{n-1}$
\begin{lem}  \label{growthygraph}
The first difference of the sequence $|\gamma_n|$ is bounded by  
$n^{p-1}$.
\end{lem}
{\bf Proof.} The edges of $\gamma_{n+1}$ are obtained from the edges
of $\gamma_n$
by applying the inverse  map $T^{-1}$. Each time a singularity  
half-line of
$T^{-1}$ intersects an edge of $\gamma_n$, this edge  splits into two,
and thus contributes $1$ to $|\gamma_{n+1}|-|\gamma_n|$.

Let $L_n$ be the set of straight lines obtained by  reflecting at most  
$n$ times
in the vertices of $P$ the extentions of the sides of $P$. By Lemma
\ref{growthgroup}, $|L_n| \prec n^{p-1}$. Each of these lines   
intersects a
singularity half-line of $T^{-1}$ at most once, therefore the total  
number of
intersections of the  lines in $L_n$ with the singularity half-lines
of $T^{-1}$
is bounded above by $n^{p-1}$. The edges of $\gamma_n$ belong to the  
lines from
$L_n$, therefore the total number of intersections of these edges with  
the
singularity half-lines of $T^{-1}$ is bounded above by $n^{p-1}$. (Note  
that the
number of the edges of $\gamma_n$ could be bigger.)  \qed

\medskip

We will now obtain the desired bound on complexity, i.e., we will
estimate the number of faces
of $\Gamma_n$. Denote by $|F_n|, |E_n|, |V_n|$ the number of faces,  
edges,
vertices of the graph $\Gamma_n$ respectively. By  
Lemma~\ref{growthygraph},
growth of the second difference of the sequence $|E_n|$ is at  most  
polynomial
of
degree $p-1$, hence $|E_n| \prec n^{p+1}$. The edges of $\Gamma_n$ are  
parallel
to the sides of $P$, thus may have at most  $p$ possible directions.
Therefore, each face of $\Gamma_n$ is at most a $2p$-gon, and
the valence of each vertex of $\Gamma_n$ is at most $2p$.
Thus, $|E_n|\leq p |F_n|,\,|E_n|\leq p |V_n|$. Euler's  formula
$|V_n|-|E_n|+|F_n|=0$ implies $p|F_n|\leq (p-1)|E_n|$, hence
$|F_n| \prec |E_n|$. \qed

\medskip

We have obtained our bound, assuming that
the rank of the abelian group generated by the sides of $P$ is $p-1$,  
i.e., maximal
possible. Although generically this is the case, the rank may drop.
Our argument proves, in fact, the statement below.

\begin{corol} \label{bound_cor}
Let $P$ be a convex euclidean $p$-gon, and let $r\leq p-1$ be the rank  
of
the abelian group generated by translations in the sides of $P$.
Then  the complexity of the outer billiard about $P$ is bounded by  
$n^{r+2}$.
\end{corol}

A polygon is {\em rational} if the rank above is $2$.
Rational polygons are dense in  the space of all polygons.
We will study
complexity of the outer billiard about a  rational polygon.
First, we recall preliminaries.

\begin{figure}[htbp]
\begin{center}
\input{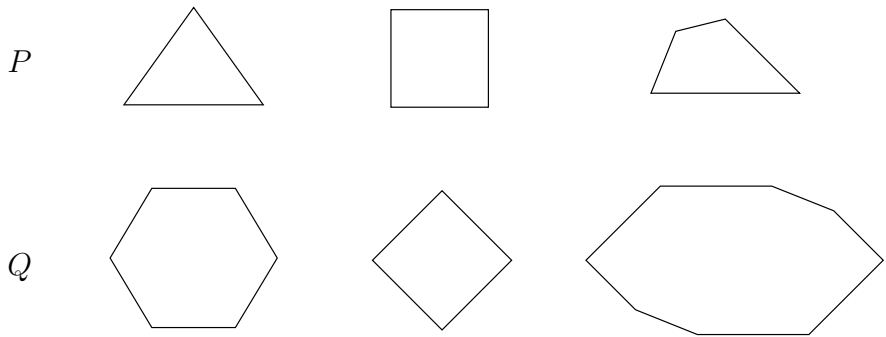}
\caption{Outer billiard; examples of polygons $P$ and $Q$}
\label{figA}
\end{center}
\end{figure}

We regard the plane as a vector space, with the center
in the interior of the convex $p$-gon $P$.
A well known construction \cite{Ta95} associates with $P$
a  homothetic family of centrally symmetric convex polygons with at most
(resp. exactly) $2p$ sides (resp. if $P$ is a generic $p$-gon).
Let $Q$ be a particular polygon in this family.
Each of its sides is parallel to a diagonal of $P$.
See figure \ref{figA}.
We endow the plane with a Minkowski norm such that $Q$ is the unit disc.
The vector norm $|\cdot|$, radius, etc, will be understood with respect
to it. We set $Q(r)=r\cdot Q$.

\begin{figure}[htbp]
\begin{center}
\input{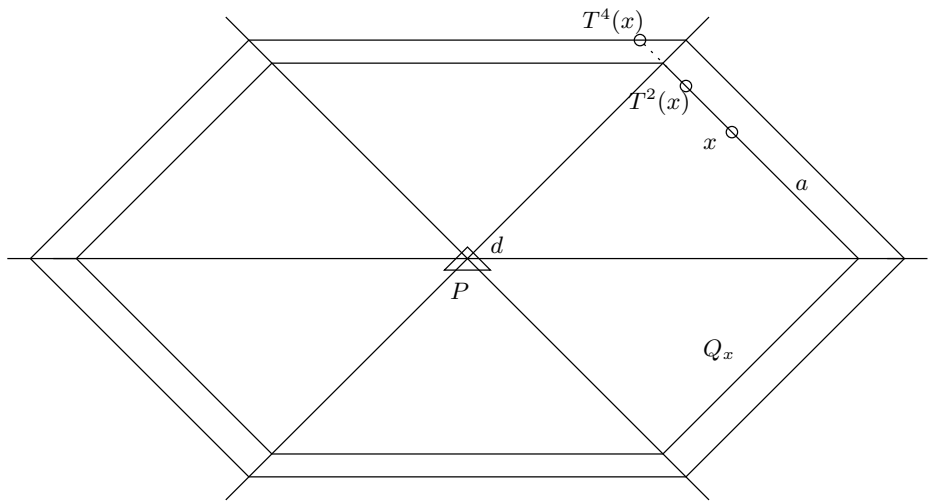}
\caption{Second iteration of the outer billiard map ``at infinity"}
\label{figB}
\end{center}
\end{figure}

The polygon $Q$ determines the geometry of orbits of $T^2$ ``at
infinity" \cite{Ta95}. We will elaborate.
Let $x$ be a point in the plane which is sufficiently far from the  
origin.
Let $Q_x$ be the circle centered at the origin and passing through $x$.
Let $a\subset Q_x$ be  the side containing $x$, and let $d$ be the
corresponding diagonal of $P$ (parallel to $a$). Then $T^2$  
translates $x$ along $a$
by $2|d|$; this continues until the orbit of $x$ overshoots $a$.
Let $y=T^{2m}$ be the corresponding point. Then the recipe above is  
applied to $y$, etc. See figure \ref{figB}.

Let $a$ be an arbitrary side of $Q$, and let $d$ be the
corresponding diagonal of $P$. The polygon $P$ is
{\em quasirational} if, up to a common  factor, the $p$ numbers  
$|a|/|d|$ are
rational.

\begin{thm}  \label{rational}
Let $P$ be a rational polygon, and let $f(\cdot)$
be the complexity of the outer billiard about $P$. Then $f(n) \sim n^2$.
\end{thm}
{\bf Proof.} Every rational polygon is quasirational.
By a construction of R. Kolodziej \cite{Ko},
there is a nested sequence of $T$-invariant polygonal simply   
connected
domains $\cdots
\subset U_i\subset U_{i+1} \subset \cdots$ exhausting the  
plane.\footnote{
This construction was used in \cite{Ko} to prove the boundedness of all  
outer
billiard orbits about quasirational  polygons.
If this is the case for arbitrary $p$-gons remains an open question for  
$p\ge 4$
\cite{Ta95}.}
By \cite{Ko}, there exists a constant $C=C(P)>0$ such that the
{\em Kolodziej domains} satisfy $Q(Ci)\subset U_i\subset Q(C(i+1))$.
We will need a general lemma.
\begin{lem} \label{infinity}
Let $P$ be an arbitrary convex polygon, and let $f(\cdot)$ be the  
complexity
of the outer billiard about $P$.
There exists $C_1>0$ such that the
contribution to $f(n)$ of the exterior of the disc of radius $C_1 n$  
grows
linearly.
\end{lem}
{\bf Proof.} We will use the preceding notation and terminology.
For any side $a$ of $Q$ set
$r_a=|a|/|d|$. Let $C_1>2/r_a$ for all sides of $P$.

Consider the $T^2$-orbit of length $n$ of an arbitrary point $x$
outside of $Q(C_1 n)$. It  follows a side, $a$, of $Q$
for $k\le n$ iterations, then it
``jumps" to the adjacent  side, $a^{\prime}$, and follows it for $n-k$
iterations.
Counting the possibilities (and assuming that $Q$ is a $2p$-gon,
i.e., that we are in the generic sitiation) we obtain $2p(n+1)$ types of
$T^2$-orbits  of length $n$. But different types mean different  
contributions
to
$f(2n)$, and vice versa.\qed

\medskip

Since $P$ is a rational polygon, the group $G\subset Iso(\RR)$
is discrete. The  graphs $\Gamma_n$ are obtained from a finite  
collection
of half-lines by $G$-action, hence $\Gamma_{\infty}=\cup_{n\ge  
1}\Gamma_n$
belongs
to a discrete collection of lines. Therefore $\Gamma_{\infty}$ is a  
graph,
and the sequence $\Gamma_1\subset\dots\subset\Gamma_n\subset\dots$
stabilizes on compacta. Moreover, there is a finite collection of convex
polygons, such that every face of $\Gamma_{\infty}$ is congruent to a polygon
in this collection. Hence the areas of the faces of $\Gamma_{\infty}$
are bounded away from zero and infinity.

Note that the constant $C_1$ in
Lemma~\ref{infinity} can be chosen arbitrarily large. We choose it so  
that
$\frac{C_1}{C}=\tau\in\N$. Then for all $n$ sufficiently large
\begin{equation}  \label{inclu_eq}
Q(C_1n)\subset U_{\tau n}\subset Q(C_1n+C).
\end{equation}

By Lemma~\ref{infinity}, up to a linear term,
$f(n)$ is the number of faces of $\Gamma_n$ intersecting $Q(C_1n)$.
By the left inclusion in (\ref{inclu_eq}), this is less than or equal
to the number of faces of $\Gamma_{\infty}$ in $U_{\tau n}$.
By preceding remarks, there is $C_2>0$ such that that number is bounded
by $C_2\ar(U_{\tau n})$. By the right inclusion in (\ref{inclu_eq}),
$\ar(U_{\tau n})$ is quadratic in $n$.
We have obtained the upper bound $f(n)\prec n^2$.

Now for the lower bound. All regular points in $X$ are periodic
\cite{Ko,GS92}. A face $F\subset X$ of $\Gamma_k$ is {\em stable}
if $F$ is a face of $\Gamma_{\infty}$.
Let $V_n\subset X$ be
the set of points  with period at most $n$. Each connected component
of $V_n$ is an open, stable face of $\Gamma_n$. By remarks above,
the number of connected components of $V_n$ has the same growth
as the area of $V_n$, thus $\ar(V_n)\prec f(n)$.
By Proposition~\ref{period_lem} below, $\ar(V_n)\sim n^2$.  \qed

\medskip

The following proposition is used in the proof of  
Theorem~\ref{rational}.
It is also of independent interest. If $g,h$ are positive functions on
$Y\subset\RR$,
the notation $g\prec h$  means that $g(x)/h(x)$
is bounded as $|x|\to\infty $. The notation  $g\sim h$ means that
$g\prec h,h \prec g$.
\begin{prop} \label{period_lem}
Let $P$ be a convex polygon and let $X_{\rm{per}}\subset X$
be the set of periodic points of the outer billiard.
For $x\in X_{\rm{per}}$ let $p(x)$ be the period.

\noindent 1. We have $|x| \prec p(x)$.

\noindent 2. Let $P$ be a rational polygon. Then for all regular points
$p(x) \sim |x|$. Let $V_n$ be the set of points such that $p(x)\le n$.
Then $\ar(V_n)\sim n^2$.
\end{prop}
{\bf Proof.} We assume without loss of generality that
$p(x)=2m$.     Let $Q_x$ be the circle through point $x$.
The sequence $x,T^2(x),\dots$ roughly follows $Q_x$.
To come back to $x$, the sequence has to go around $Q_x$ at least once.
Let $\delta$ be the ``largest step" of $T^2$. Then we need at least
$\per(Q_x)/\delta$ steps to return. Since $\per(Q_x)\sim |x|$,
the first claim follows.

Let now $P$ be rational, and hence quasirational polygon,
and let $U_k,\,k\ge 1,$ be the Kolodziej domains. Let $k=k(x)$
be such that $x\in U_k\setminus U_{k-1}$. The relations
$Q(Ck)\subset U_k\subset Q(C(k+1))$ imply that the function
$k(x)$ satisfies $k(x)\sim |x|$. By inclusion
$U_k\setminus U_{k-1}\subset Q(C(k+1))\setminus Q(C(k-1))$,
we have $\ar(U_k\setminus U_{k-1})\sim  |x|$. The point $x$
belongs to a unique face, $F=F(x)$, of $\Gamma_{\infty}$,
hence $p(x)\ar(F) \le \ar(U_k\setminus U_{k-1})$. By preceding remarks,
$p(x)\prec \ar(U_k\setminus U_{k-1})$, implying
$p(x)\prec |x|$, and hence the equivalence $p(x)\sim |x|$.

By this relation, there are constants $C_3,C_4>0$ such that for $n$  
sufficiently
large
$Q(C_3n)\subset V_n \subset Q(C_4n)$, proving the last claim.\qed

\subsection{The elliptic and the hyperbolic cases}    
\label{neg_curv_sub}
We will first study the outer billiard in elliptic geometry.

\begin{thm} \label{spheresubexp1}
Let $P\subset S^2$ be a convex spherical polygon.
The  complexity of  outer
billiards about $P$ grows subexponentially.
\end{thm}
{\bf Proof.} For $x\in S^2$ let $l=x^*$ be the appropriately oriented
great circle centered at $x$. This diffeomorphism $S^2\to Geo(S^2)$
is the {\em spherical duality}, and we denote by $x=l^*$ the inverse
diffeomorphism.

Let $A,B,\dots$ be the vertices of $P$. The geodesics  
$a=A^*,b=B^*,\dots$
bound the convex polygon $P^*$. The correspondence $P\mapsto P^*$
is an automorphism of the space of convex spherical polygons.
The proof of the following lemma is contained in \cite{Ta95a}.
(See also \cite{Ta95}.)
\begin{lem}  \label{duality}
Let $P,P^*\subset S^2$ be as above. Let $X_o,X_b$ be the phase
spaces of the outer billiard about $P$, inner billiard about
$P^*$, and let $T_o:X_o\to X_o,\, T_b:X_b\to X_b$ be the
outer billiard, inner billiard maps respectively.

The spherical duality induces a diffeomorphism $X_o \to X_b$;
it conjugates $T_o:X_o\to X_o$ and $T_b:X_b\to X_b$;
it induces an isomorphism of the coding of $T_o$-orbits by corners of  
$P$ and
the
coding of $T_b$-orbits by sides of $P^*$.
\end{lem}

\begin{figure}[htbp]
\begin{center}
\input{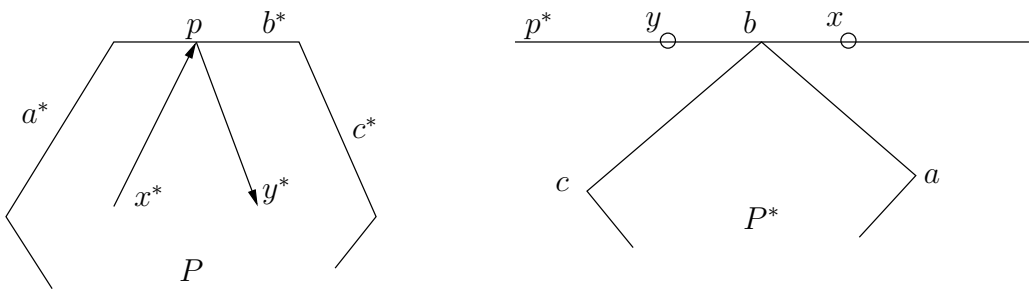}
\caption{Duality between inner and outer billiards}
\label{figH}
\end{center}
\end{figure}

Figure~\ref{figH} illustrates Lemma~\ref{duality}.
Let $f_o(n)$ (resp.  $f_b(n)$) be the corner complexity of the
outer billiard about $P$ (resp.  billiard in $P^*$).
By Lemma~\ref{duality}, $f_o(n)=f_b(n+1)$. The claim now follows from
Theorem~\ref{spheresubexp}.  \qed

\vspace{5mm}

\begin{figure}[htbp]
\begin{center}
\input{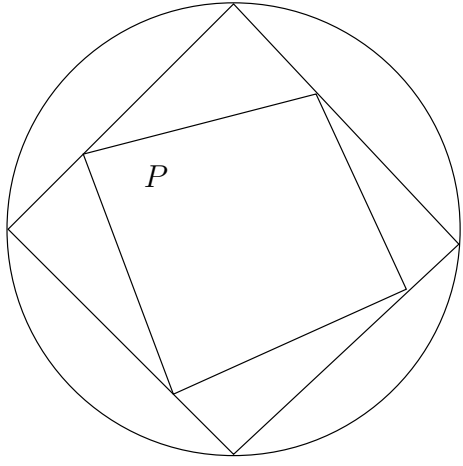}
\caption{A large quadrilateral}
\label{figC}
\end{center}
\end{figure}

We refer to \cite{DT03} for the background on polygonal outer
billiards in the hyperbolic plane. Let $P\subset\HT$ be a
$p$-gon, let $X=\HT\setminus P$, and let $T:X\to X$ be the
outer billiard map. It extends to a homeomorphism
of the circle at infinity, $\tau:S\to S$. Its rotation
number satisfies $\rho(P) \geq 1/p$ \cite{DT03}. The polygon $P$ is  
{\em large}
if $\rho(P) = 1/p$ and $\tau$ has a hyperbolic $p$-periodic orbit.
See figure \ref{figC}.
The set of large polygons is open in the natural topology \cite{DT03}.

\begin{thm}    \label{hyplowbd}
Let $P\subset\HT$ be an arbitrary convex polygon, and let $f(\cdot)$ be  
the
complexity of the outer billiard about $P$.
Then  $n \prec f(n)$. If $P$ is a large polygon, then $f(n)\sim n$.
\end{thm}
{\bf Proof.} The bound $n \prec f(n)$ fails iff the sequence  
$\Gamma_k,\,k\ge
1,$
stabilizes. Assume this to be the case, and let
$\Gamma_m=\Gamma_{m+1}=\dots=\Gamma_{\infty}$.
The outer billiard map $T:X\to X$ preserves $\Gamma_{\infty}$; the  
restriction
of $T$
to a closed face of $\Gamma_{\infty}$ is a diffeomorphism onto another  
one.
Since $\Gamma_{\infty}$ is a finite graph, we find $n\in\N$ such that every  
face
of $\Gamma_{\infty}$ is invariant under $T^n$.

\begin{figure}[htbp]
\begin{center}
\input{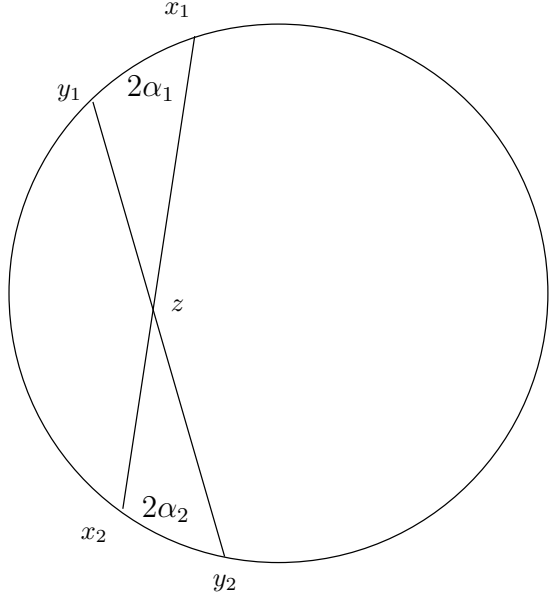}
\caption{Computing the distortion of the map $\tau$}
\label{figE}
\end{center}
\end{figure}

Let $F$ be a closed face of $\Gamma_{\infty}$. Then $\bo F\cap S$ is  
either
empty, or a vertex, or an edge of $F$. We will study the latter.
Let $v_1,\dots,v_N\in S$ be the consecutive endpoints of these edges,
let $e_i\subset S$ (resp. $\alpha_i\subset \HT$)
be the circular arc (resp. the geodesic) with endpoints $v_i,v_{i+1}$  
(we set
$N+1=1$),
and  let $F_i$ be the corresponding face of $\Gamma_{\infty}$.
The restriction $T^n|_{F_i}$ is induced by an isometry, $g_i\in  
Iso(\HT)$.
The elements $g_1,\dots,g_N$ are all equal to the identity iff  
$\tau^N=1$.
\begin{lem} \label{nonperiod}
The map $\tau:S\to S $ is not periodic.
\end{lem}
{\bf Proof.} Let $z$
be a vertex of $P$. For close points $x_1,y_1 \in S$
let $x_2,y_2 \in S$ be their reflections about $z$. Let
$\lambda_1=|x_2z|/|x_1 z|$ and let $2\alpha_i$ be the angular measure  
of the arc
$x_i y_i,\,i=1,2$; see figure \ref{figE}. The
triangles $x_1 z y_1$ and $x_2 z y_2$ are similar, therefore
\begin{equation} \label{cyclicorb}
\sin \alpha_2 = \lambda_1 \sin \alpha_1.
\end{equation}

Let $x_1,\dots ,x_N$ be a periodic trajectory of the map $\tau$  
consisting of
smooth
points, and let $\lambda_1\dots ,\lambda_N$ be the
respective ratios. Set $\Lambda=\Pi_{i=1}^N \lambda_i$.
Let $y_1$ be a point sufficiently close to $x_1$,  and let
$y_1,\dots ,y_N$ be its $\tau$-orbit; we assume that for both orbits
the reflections occur in the same vertices  of $P$. It follows from
equation~\ref{cyclicorb} that $y_1,\dots ,y_N$ is
a periodic trajectory iff $\Lambda=1.$ In particular, if $\tau$ has a  
periodic
interval, then $\Lambda=1$ there.

\begin{figure}[htbp]
\begin{center}
\input{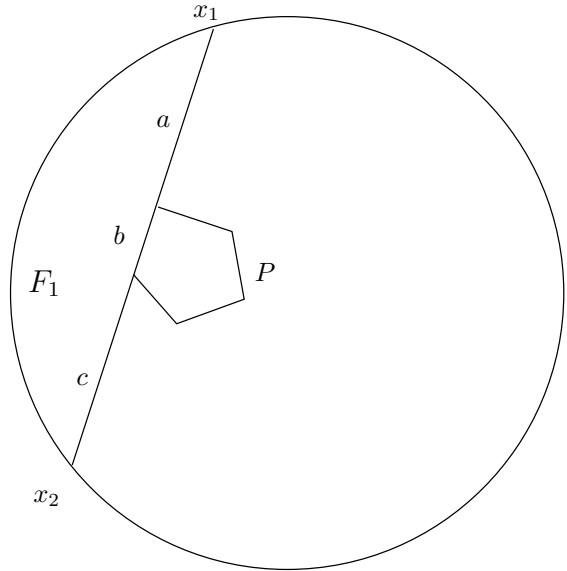}
\caption{Destruction of a periodic orbit of $\tau$}
\label{figF}
\end{center}
\end{figure}

Let now $x_1$ cross counter-clockwise a singularity
half-line of $T$. In the notation of
figure \ref{figF}, $\lambda_1=(b+c)/a$ (resp. $\lambda_1=c/(a+b)$)
right before (resp.  after) this.  By $(b+c)/a > c/(a+b)$,
the equality $\Lambda=1$
before a singularity half-line  implies that $\Lambda<1$ immediately  
after
it. \qed

\medskip

By Lemma~\ref{nonperiod}, we can assume without loss of generality 
that $g_1\ne 1$.
Then $g_1$ is a (hyperbolic) parallel translation with the axis
$\alpha_1$, and $F_1$ is the domain bounded by $\alpha_1$ and $e_1$.
We will say that $F_1$ is a {\em lunar face} of $\Gamma_{\infty}$.
The union of  lunar faces of $\Gamma_{\infty}$ is invariant under
$T$. Therefore for any $k>0$ there is $l=l(k)$ such that
$T^{-k}(\alpha_1)=\alpha_l$. A geodesic $\alpha_i,\,1\le i\le N,$  
cannot contain
a side of $P$. If it does, then $F_i$ contains a singular line of
$T$ in its interior, contrary to the definition
of $F_i$. See figure \ref{figF} where $x_1 x_2$ represents now the geodesic $\alpha_1$.
Thus, $\alpha_1$ is not an edge of $\Gamma_m$ for any $m$.
This contradiction proves our first claim.

\begin{figure}[htbp]
\begin{center}
\input{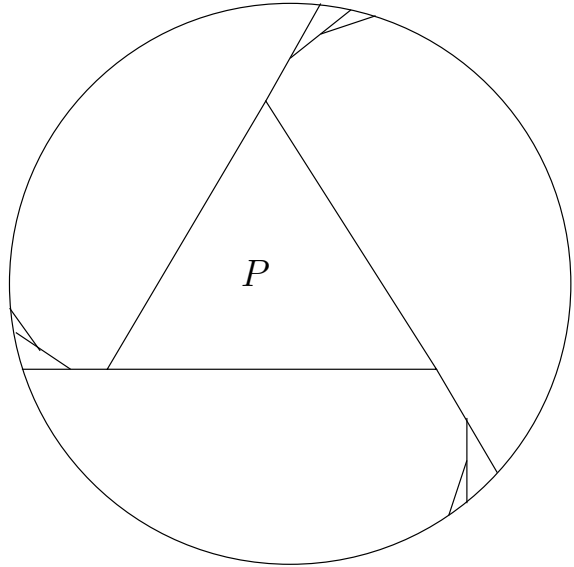}
\caption{The graph $\Gamma_2$ for a large triangle}
\label{figD}
\end{center}
\end{figure}

Let now $P$ be a large $p$-gon.
Then $\Gamma_n$ is a disjoint union of $p$ binary
trees \cite{DT03}  (see figure~\ref{figD}), and
hence $|\Gamma_n|$ grows linearly.  \qed

\begin{figure}[htbp]
\begin{center}
\input{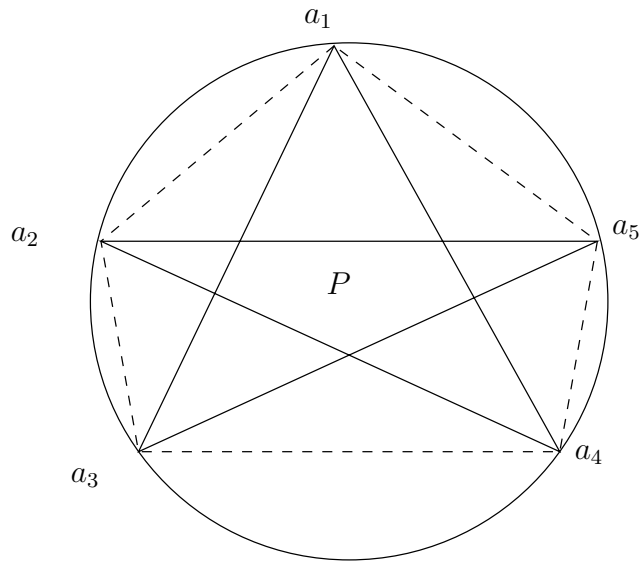}
\caption{Finite complexity of the outer billiard map at infinity}
\label{figG}
\end{center}
\end{figure}

\begin{rem} \label{degen_rem}
{\rm   Note that the function $f(\cdot)$
is bounded below by the complexity of the induced
map $\tau:S\to S$ with respect to the natural partition.
However, the latter may be finite. See figure~\ref{figG}. }
\end{rem}

\end{document}